# TRANSITION FROM THE ANNEALED TO THE QUENCHED ASYMPTOTICS FOR A RANDOM WALK ON RANDOM OBSTACLES

By Gérard Ben Arous, Stanislav Molchanov and Alejandro F. Ramírez[1]

*New York University, University of North Carolina–Charlotte and Pontificia Universidad Católica de Chile*

In this work we study a natural transition mechanism describing the passage from a quenched (almost sure) regime to an annealed (in average) one, for a symmetric simple random walk on random obstacles on sites having an identical and independent law. The transition mechanism we study was first proposed in the context of sums of identical independent random exponents by Ben Arous, Bogachev and Molchanov in [*Probab. Theory Related Fields* **132** (2005) 579–612]. Let $p(x,t)$ be the survival probability at time $t$ of the random walk, starting from site $x$, and let $L(t)$ be some increasing function of time. We show that the empirical average of $p(x,t)$ over a box of side $L(t)$ has different asymptotic behaviors depending on $L(t)$. T here are constants $0<\gamma_1<\gamma_2$ such that if $L(t)\geq e^{\gamma t^{d/(d+2)}}$, with $\gamma>\gamma_1$, a law of large numbers is satisfied and the empirical survival probability decreases like the annealed one; if $L(t)\geq e^{\gamma t^{d/(d+2)}}$, with $\gamma>\gamma_2$, also a central limit theorem is satisfied. If $L(t)\ll t$, the averaged survival probability decreases like the quenched survival probability. If $t\ll L(t)$ and $\log L(t)\ll t^{d/(d+2)}$ we obtain an intermediate regime. Furthermore, when the dimension $d=1$ it is possible to describe the fluctuations of the averaged survival probability when $L(t)=e^{\gamma t^{d/(d+2)}}$ with $\gamma<\gamma_2$: it is shown that they are infinitely divisible laws with a Lévy spectral function which explodes when $x\to 0$ as stable laws of characteristic exponent $\alpha<2$. These results show that the quenched and annealed survival probabilities correspond to a low-

Received June 2003; revised November 2004.
[1]Supported in part by Fondo Nacional de Desarrollo Científico y Tecnológico Grant 1020686.

and high-temperature behavior of a mean-field type phase transition mechanism.

**1. Introduction.** When studying the long-time behavior of Markovian dynamics in random media, one has to choose between two possibilities. One can either choose to fix an almost sure realization of the random medium in which the dynamics takes place, and then let the time go to infinity, or first average the randomness of the medium before letting the time grow. We will call the first approach the *quenched* regime, and the second one the *annealed* or more appropriately the *averaged* regime. It is often the case that the two approaches give completely different answers. To name a few instances of such problems let us mention the problem of random walks among random traps (or Brownian motion among Poissonian traps) (see [1, 4] or Chapter 4 of [12]), or slowdowns of random walks in random environments [9, 10], or the phenomenon of intermittency for the parabolic Anderson model [6, 7].

We want to address here the question of the relevance of these two approaches. A reasonable first answer to this question could simply be that the quenched approach solves the true question but that the averaged approach, being often much simpler, has the merit of being the first possibility to understand a hard problem (historically this has been the case with all the examples quoted above). We will here introduce another way to address this question by showing that there is a rich transition between these two approaches which shows that they are the two extreme points of a broad range.

To introduce this new transition between the quenched and annealed regimes we need to introduce, at least, one new parameter. This parameter will be a spatial scale, that is, the spatial extent, say $L$, of the initial data. We will consider a quenched realization of the medium and will let this spatial scale $L$ and the time $t$ tend to infinity together. It is rather clear that if the diverging time scale is small enough as a function of $L$, or equivalently if the spatial scale $L$ is large enough as a function of $t$, the annealed regime should prevail. Indeed the random medium should be thought as sampled enough by the initial data for a spatial ergodic theorem to apply and thus justify the annealed asymptotics. On the other hand, on much longer time scales, when the time is very large to make the spatial scale look very small and not very different from a pointwise initial datum, the quenched asymptotics should be in force. We will show that there is a rich transition of asymptotics regimes interpolating between these two extremes. This transition should be seen as a transition between the *bulk* properties and the extreme values of certain local spectral statistics of the random medium. One should thus use the results provided by the annealed or quenched asymptotics with care depending on this initial spatial scale.

This point of view was first proposed and studied in [3] by Ben Arous, Bogachev and Molchanov for the very simple case of large sums of random exponentials. We have chosen to work here in the context of random walks among random traps. We will address the question of branching random walks and slowdowns of random walks in random environments in other works. In this article we study the transition mechanism proposed in [3] between the quenched and annealed behavior of a random walk on random traps on the lattice.

More specifically, we study the asymptotic behavior of the survival probability of a random walk which is killed when touching traps distributed according to a product Bernoulli random variable. This can be regarded as the case of the parabolic Anderson problem where the random potential takes only two values: 0 or $-\infty$.

Let us call $p(x,t,w)$ the probability that a simple symmetric continuous-time random walk survives up to time $t$, starting from site $x$, on a random trap environment given by the configuration $w=\{w_x: x\in\mathbf{Z}^d\}\in\{0,1\}^{\mathbf{Z}^d}$ (here if $w_x=1$ at a site $x$, then $x$ corresponds to a hard trap, killing the random walk the first time it hits $x$). In [4, 5] Donsker and Varadhan showed that for long times, the average of $p(0,t,w)$ with respect to a Bernoulli product measure of parameter $p$ behaves like $\exp\{-c_2(d,p)t^{d/(d+2)}+o(t^{d/(d+2)})\}$, where $c_2(d,p)$ is a constant depending only on the dimension and $p$. This is the annealed behavior. Subsequently, in [2] Antal, following Sznitman, showed that almost surely with respect to the distribution of the traps, and on the event that the origin belongs to an infinite trap free component, $p(0,t,w)$ behaves like $\exp\{-c_1(d,p)t/(\log t)^{d/2}+(t/(\log t)^{d/2})\}$, where $c_1(d,p)$ is another constant depending on the dimension and the parameter $p$. This is the quenched behavior. In this paper we study the averaged quantity $p^L(0,t,w):=\frac{1}{|\Lambda_L|}\sum_{x\in\Lambda_L}p(x,t,w)$, where $\Lambda_L:=[-(2L+1),(2L+1)]^d\cap\mathbf{Z}^d$ is a box of radius $L>0$ and $|\Lambda_L|$ is its cardinality. If we make $L$ depend on $t$, then as $t\to\infty$ we show that several behaviors occur, depending on the rate at which $L$ grows with $t$. This is the content of Theorems 1 and 2 stated in the following section. In particular, it is shown that for $L\le t$, the averaged survival probability $p^L(0,t,w)$ behaves as in the quenched situation, while if $L\ge\exp\{\gamma\frac{c_2}{d}t^{d/(d+2)}\}$ it behaves as in the annealed case. In the one-dimensional case it is possible to give more precise results for $L=\exp\{\gamma\frac{c_2}{d}t^{d/(d+2)}\}$. This is the content of Theorem 3.

Let us now briefly discuss the intuitive picture described behind the proof of Theorem 1. It is shown there that there exist four main averaging regimes (see also Table 1): case 1 or $L(t)\le t$; case 2 or $L(t)\ge t$ and $\log L(t)\ll t^{d/(d+2)}$, where for two real functions defined on $[0,\infty)$, $f\ll g$ means that $\lim_{t\to\infty}f(t)/g(t)=0$; case 5 or $L(t)\ge e^{(\gamma/d)c_2t^{(d/d+2)}}$, with $\gamma>2/(d+2)$ and $c_2$ a constant depending on $d$ and $p$; and case 6, or $L(t)\ge e^{(\gamma/d)c_2t^{(d/d+2)}}$, with $\gamma > 2^{d/(d+2)} 2/(d+2)$. In cases 1 and 2, basically, in the average defining the survival probability $p^L$ there is a dominant clearing in the random environment giving the main contribution to the logarithmic asymptotic behavior. This is the clearing of radius proportional to $(\log L)^{1/d}$. Thus, the difference between cases 1 and 2 corresponds simply to the case in which $\log L$ is of the order of $\log t$ (case 1) and $\log L \gg \log t$ (in general in case 2). On the other extreme, cases 5 and 6 are situations where the behavior of $p^L$ is determined by many large clearings of the random environment, it being impossible to isolate a single dominant one.

In Section 2 of this paper, the results are stated. In Section 3 several important estimates, concerning the survival probabilities analyzed as a field of random variables, are derived. Most of them are moment and correlation estimates. In Section 4, cases 1 and 2 of Theorem 1 are derived, using the method of enlargement of obstacles of Sznitman (Chapter 4 of [12]). In Section 5 the cases 3, 4, 5 and 6 are derived. The general philosophy of the corresponding proofs is the use of renormalization methods to control the main contributions of the random environment. The proof of Theorem 3 is the content of Section 6. The tools used here correspond to standard, though lengthy, verifications of the necessary hypothesis for the convergence of a given sequence of random variables to an infinitely divisible law.

## 2. Notation and results.

We will first introduce the necessary notation to define a symmetric simple random walk on the lattice $\mathbb{Z}^d$ of total jump rate 1 in a random obstacle environment. We define a random environment through a product measure $\mu$ on the Cartesian product $X := \{0,1\}^{\mathbb{Z}^d}$ with the Borel-$\sigma$ field generated by the product topology, so that $\mu(w(x) = 1) = p$, where $0 < p < 1$ and $w(x)$ is the $x$ coordinate of $w \in X$. Each element of $w \in X$ will be called an *obstacle environment*. A site $x$ of the lattice where $w(x) = 1$ represents a site with an obstacle, while if $w(x) = 0$ there is none. Given any real function $f(w)$ of the environment $w$ we will, throughout this paper, denote as $\langle f \rangle := \int f(w)\, d\mu$ the expectation of $f$ with respect to the law of the environment. Let us now denote by $\mathcal{G}(w) := \{y \in \mathbb{Z}^d : w(y) = 1\}$ the set of sites having an obstacle or *obstacle set*. Throughout the sequel $Z_{\cdot}$ will denote the canonical $d$-dimensional continuous-time symmetric simple random walk of total jump rate 1 defined on the Skorokhod space $D([0,\infty); \mathbb{Z}^d)$. We will call $P_x$ the law of such a random walk starting from site $x \in \mathbb{Z}^d$ and $\tau(w) := \inf\{t \geq 0 : Z_t \in \mathcal{G}(w)\}$ the killing time, or the first hitting time of the obstacle set $\mathcal{G}(w)$. Let us now denote as $p(x,t,w) := P_x(\tau(w) > t)$ the probability that a random walk starting from site $x$ does not hit the obstacle set $\mathcal{G}(w)$ by time $t$. Such a probability will be referred to as the *quenched survival probability* at time $t$ of a random walk starting from site $x$. Similarly we will call $\langle p(x,t) \rangle := \int p(x,t,w)\, d\mu$ the *annealed survival probability* at time $t$ of a random walk starting from site $x$. Furthermore, we denote the sets $\{p(x,t,w) : x \in \mathbb{Z}^d\}$ and $\{\langle p(x,t)\rangle : x \in \mathbb{Z}^d\}$ as the *field of quenched survival probabilities* and the *field of annealed survival probabilities*, respectively. In the sequel, whenever there is no danger of confusion, we will drop the variables $x$, $w$ and $t$ of the survival probabilities, writing $p(x,t)$ or simply $p$ in place of the quenched $p(x,t,w)$ and writing $\langle p\rangle$ in place of the annealed $\langle p(x,t)\rangle$. Given $r \in [0,\infty)$ and $x \in \mathbb{Z}^d$ we denote by $\Lambda(x,r) := \{y \in \mathbb{Z}^d : \|x-y\| \le r\}$ the ball of radius $r$ centered at site $x$ under the norm $\|x\| := \sup_{i=1,\dots,d} |x_i|$, where $x_i$ are the coordinates of site $x$. We will frequently write $\Lambda_r$ in place of $\Lambda(0,r)$. In contrast we will denote by $B(x,r) := \{y \in \mathbb{Z}^d : |x-y| \le r\}$ the Euclidean ball of radius $r$ centered at site $x$ under the Euclidean norm $|x| := \sqrt{\sum_{i=1}^d x_i^2}$.

Now that we have settled the basic notation, we are in a position of defining quantities which will correspond to a transition between the quenched and annealed survival probabilities. We define the *averaged survival probability* at scale $L$ and time $t$ for a random walk starting from site $x$ as

$$p^L(x,t,w) := \frac{1}{|\Lambda(x,L)|} \sum_{y \in \Lambda(x,L)} p(y,t,w),$$

where $|U|$ denotes the cardinality of $U \subset \mathbb{Z}^d$. Whenever there is no danger of confusion, we will drop the variables $x$, $t$ or $w$, writing $p^L$ or $p^L(x,t)$ in place of the averaged survival probability $p^L(x,t,w)$.

Given two real-valued functions $f, g$, the notation $f \sim g$ means that $\lim_{u\to\infty} \frac{g(u)}{f(u)} = 1$, $f \gg g$ means that $\lim_{u\to\infty} \frac{g}{f} = 0$, while $f \ll g$ means that $g \gg f$. Also, $w_d$ and $\ell_d$ will denote the volume of the ball of unit radius in $\mathbb{R}^d$ and the principal Dirichlet eigenvalue of the differential operator $-\frac{1}{2d}\Delta$ on this ball, respectively. For $p<1$ we define $\nu := |\log(1-p)|$ and the constants $c_1(d,p) := \ell_d/R_0^2$ with $R_0 := (\frac{d}{w_d \nu})^{1/d}$ and $c_2(d,p) := (w_d\nu)^{2/(d+2)}((d+2)/2)(2\ell_d/d)^{d/(d+2)}$. Next, we define $p_c := \inf\{p : \mu(|\mathcal{G}(w)| = \infty) > 0\}$ the critical probability of site percolation for the obstacles on $\mathbf{Z}^d$.

In this paper we want to study the behavior of the averaged survival probabilities $p^L(0,t,w)$ for large $t$ and $L$. By standard ergodic theorems it is possible to show that for $t$ fixed, as $L \to \infty$ we have that $\mu$-a.s. $p^L(0,t,w) \sim \langle p\rangle$ and the behavior is annealed. On the other hand, it is not difficult to see that for $L$ fixed, as $t \to \infty$ we have that $\mu$-a.s. $p^L(0,t,w) \sim p(0,t,w)$ and the behavior is quenched. It is natural to ask if there is a transition mechanism between these two extremes when we let both $L \to \infty$ and $t \to \infty$. In this paper we partially answer this question. We let $L = L(t)$ depend on time $t$ so that $L(t) \gg 1$ and distinguish six cases according to the growth rate of $L(t)$: If $L(t) \le t$ we say the asymptotics is in *case* 1; if $\log L(t) \ll t^{d/(d+2)}$ and $L(t) \ge t$ we say the asymptotics is in *case* 2; if $L(t) = \exp\{\gamma c_2 t^{d/(d+2)}/d\}$ with $\gamma < \gamma_1 := 2/(d+2)$ we say the asymptotics is in *case* 3; if $L(t) = \exp\{\gamma c_2 t^{d/(d+2)}/d\}$ with $\gamma_1 < \gamma < \gamma_2 := 2^{d/(d+2)}\gamma_1$ we say the asymptotics is in *case* 4; if $L(t) \geq \exp\{\gamma c_2 t^{d/(d+2)}/d\}$ with $\gamma > \gamma_1$ we say the asymptotics is in *case* 5; while if $L(t) \geq \exp\{\gamma c_2 t^{d/(d+2)}/d\}$ with $\gamma > \gamma_2$ we say the asymptotics is in *case* 6. We summarize this classification in Table 1.

We now state the main result of this paper.

THEOREM 1. *Let $L(t) : [0, \infty) \to \mathbf{N}$ be some nondecreasing function and assume that $0 \leq p < 1$. Then the following statements are true.*

(i) Case 1. *If $1 \ll L(t) \leq t$ and $p < 1 - p_c$, then $\mu$-a.s. we have that*

$$\log p^{L(t)}(0, t, w) \sim -c_1(d, p)\frac{t}{(\log t)^{2/d}}.$$

(i′) Case 2. *If $L(t) \geq t$ and $\log L(t) \ll t^{d/(d+2)}$, then $\mu$-a.s. it is true that*

$$\log p^{L(t)}(0, t, w) \sim -c_1(d, p)\frac{t}{(\log L(t))^{2/d}}.$$

(ii) Case 5. *Let $\gamma_1 = \frac{2}{d+2}$ and suppose that there is a $\gamma > 0$ such that $L(t) \geq \exp(\gamma \frac{c_2}{d} t^{d/(d+2)})$. If $\gamma > \gamma_1$, then in $\mu$-probability it is true that*

$$\frac{p^{L(t)}(0, t, w)}{\langle p(0, t)\rangle} \sim 1. \tag{1}$$

*In particular, in $\mu$-probability we have that*

$$\log p^{L(t)}(0, t, w) \sim -c_2(d, p)t^{d/(d+2)}. \tag{2}$$

*On the other hand, if $L(t) \leq \exp(\gamma \frac{c_2}{d} t^{d/(d+2)})$ with $\gamma \leq \gamma_1$, then $\mu$-a.s. it is true that*

$$\frac{p^{L(t)}(0, t, w)}{\langle p(0, t)\rangle} \ll 1. \tag{3}$$

TABLE 1

| | | |
|---|---|---|
| Case 1 | $L(t) \leq t$ | |
| Case 2 | $L(t) \geq t$ | $\log L(t) \ll t^{d/(d+2)}$ |
| Case 3 | $L(t) = e^{(\gamma/d)c_2 t^{d/(d+2)}}$ | $\gamma < \gamma_1$ |
| Case 4 | $L(t) = e^{(\gamma/d)c_2 t^{d/(d+2)}}$ | $\gamma_1 < \gamma < \gamma_2$ |
| Case 5 | $L(t) \geq e^{(\gamma/d)c_2 t^{d/(d+2)}}$ | $\gamma > \gamma_1$ |
| Case 6 | $L(t) \geq e^{(\gamma/d)c_2 t^{d/(d+2)}}$ | $\gamma > \gamma_2$ |

(iii) Case 6. *Let* $\gamma_2 = 2^{d/(d+2)}\frac{2}{d+2}$. *If* $L(t) \geq \exp(\gamma \frac{c_2}{d} t^{d/(d+2)})$ *for some* $\gamma > \gamma_2$, *then*

$$\lim_{t\to\infty} \frac{p^{L(t)}(0,t,w) - \langle p(0,t)\rangle}{\sqrt{\mathrm{Var}_\mu(p^{L(t)}(0,t,w))}} = \mathcal{N}(0,1), \tag{4}$$

*where for any random variable* $X(w)$, $\mathrm{Var}_\mu(X) = \int (X - \int X\, d\mu)^2\, d\mu$ *and* $\mathcal{N}(0,1)$ *is the normal law of variance* 1 *and mean* 0 *and the convergence is in the sense of distributions. On the other hand, if* $L(t) \leq \exp(\gamma \frac{c_2}{d} t^{d/(d+2)})$ *for some* $\gamma_1 < \gamma < \gamma_2$, *then in* $\mu$*-probability it is true that*

$$\frac{p^{L(t)}(0,t,w) - \langle p(0,t)\rangle}{\sqrt{\mathrm{Var}_\mu(p^{L(t)}(0,t,w))}} \ll 1. \tag{5}$$

It is possible to complement the statement of Theorem 1 in the case in which the scale $L(t)$ is of the form $\exp\{\gamma \frac{c_2}{d} t^{d/(d+2)}\}$. For reasons that will become clear below, we define the concept of negative conjugate constants. Let $x, y$ be positive real numbers. We say that $x$ is *negative conjugate* to $y$ if $\frac{1}{x} - \frac{1}{y} = 1$. Let us now define

$$\alpha := \frac{d}{2}, \tag{6}$$

$$\alpha' := \frac{d}{d+2}. \tag{7}$$

Clearly, $\alpha'$ is negative conjugate to $\alpha$.

Theorem 2. *Let* $\gamma > 0$ *and* $a(\gamma) := \alpha'(\frac{\alpha}{\alpha'}\gamma)^{-1/\alpha} + \gamma$ *and* $L(t) := e^{\gamma(c_2/d)t^{d/(d+2)}}$. *Then the following statements are true.*

(i) Case 3. *If* $\gamma \leq \gamma_1 = \frac{2}{d+2}$, *then for every* $0 < \delta < 1$, $\mu$*-a.s. we have that*

$$\frac{p^{L(t)}(0,t,w)}{\exp\{-(a(\gamma) - \delta)c_2 t^{d/(d+2)}\}} \ll 1. \tag{8}$$

(ii) Case 4. *If* $\gamma_1 < \gamma < \gamma_2 = 2^{d/(d+2)}\frac{2}{d+2}$, *then for every* $0 < \delta < 1$ *we have in* $\mu$*-probability that*

$$\frac{p^{L(t)}(0,t,w) - \langle p(0,t)\rangle}{\exp\{-(a(\gamma) - \delta)c_2 t^{d/(d+2)}\}} \ll 1. \tag{9}$$

We do expect the function $a(\gamma)$ to be critical in the sense that if $\gamma \leq \gamma_1$, then for every $0 < \delta < 1$, $\mu$-a.s. we should have that $\frac{p^{L(t)}(0,t,w)}{\exp\{-(a(\gamma)+\delta)c_2 t^{d/(d+2)}\}} \gg 1$. We also expect a similar statement complementing case (ii) of Theorem 2.

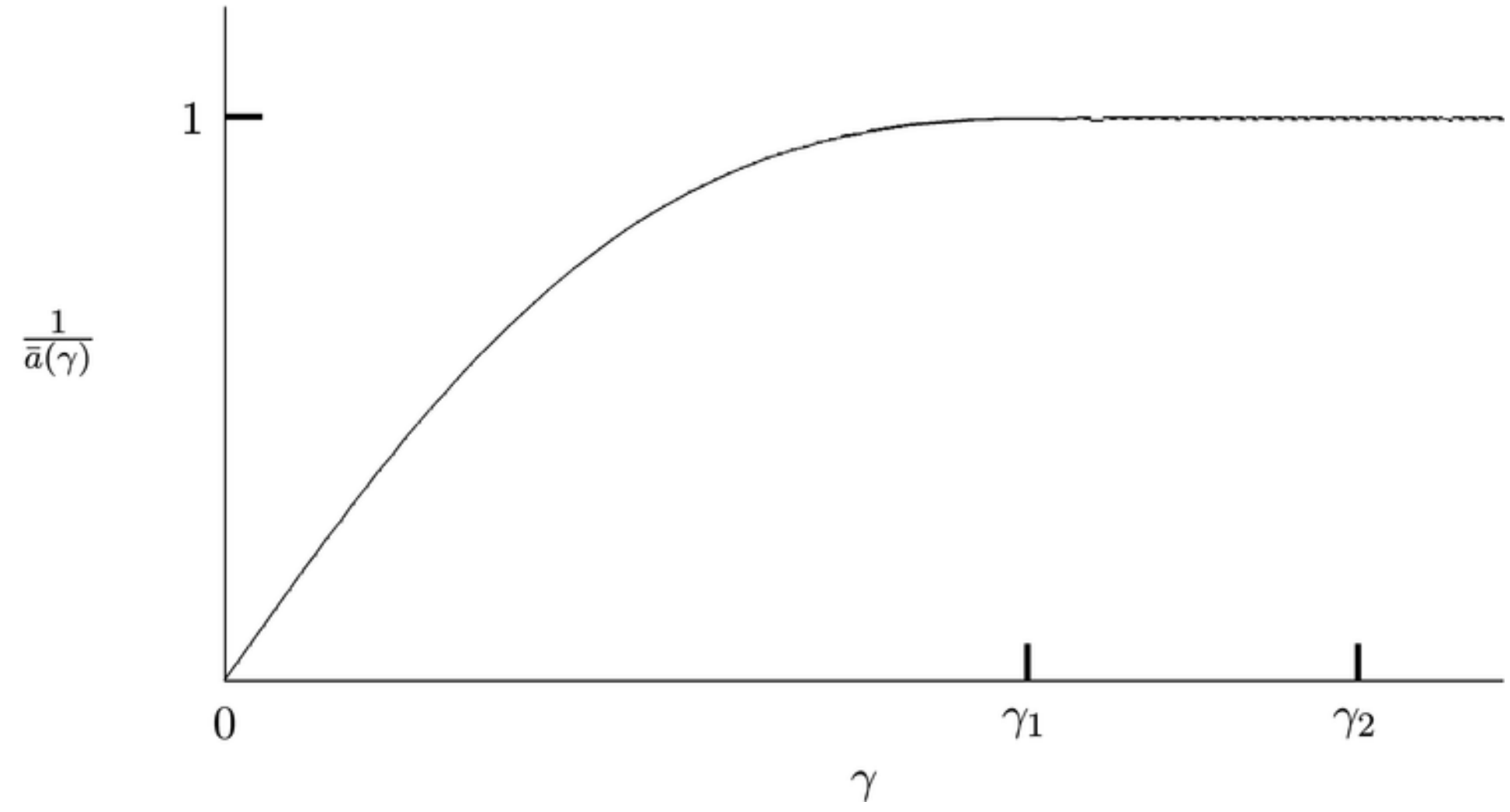

FIG. 1. *Plot of* $1/\bar{a}(\gamma)$ *when* $d=2$. *Here* $\gamma_1 = 0.5$ *and* $\gamma_2 = 0.71$.

The statements of Theorems 1 and 2 can be visualized as a plot in the case in which the time-dependent scale $L(t)$ of the averaged survival probability $p^{L(t)}(0,t,w)$ is of the form $L(t) = \exp\{\gamma \frac{c_2}{d} t^{d/(d+2)}\}$. In fact, these results suggest that the averaged survival probability $p^{L(t)}(0,t,w)$ behaves as $\log p^{L(t)}(0,t,w) \sim -\bar{a}(\gamma) t^{d/(d+2)}$ with $\bar{a}(\gamma) := c_2(d,p)$ for $\gamma > \gamma_1$, $\bar{a}(\gamma) := a(\gamma) = \frac{d}{d+2}(\frac{d+2}{2}\gamma)^{-2/d} + \gamma$ for $0 < \gamma < \gamma_1$ while $\bar{a}(\gamma) := 0$ for $\gamma = 0$. This is illustrated in Figure 1 in the case in which the dimension of the lattice has the value $d = 2$.

Figure 1 shows graphically the regions where the law of large numbers ($\gamma > \gamma_1$) and the central limit theorem ($\gamma > \gamma_2$) are valid. The quantity $a(\gamma)c_2$ plays the role of a "free energy" and the graph shows the presence of a phase transition at $\gamma = \gamma_1$, which is a point of nonanalyticity. This kind of behavior is analogous to that of mean-field statistical mechanics magnetization models such as the random energy model (REM). However, it should be pointed out that the constancy of $a(\gamma)$ for $\gamma > \gamma_1$ does not correspond to a "freezing" phenomenon as observed in the REM. In fact, the constancy of $a(\gamma)$ for $\gamma > \gamma_1$ can be interpreted as a "high-temperature" phenomenon coming from a law of large numbers, whereas the freezing phenomenon of the REM is a low-temperature phenomenon where the main contribution to the free energy comes from low energy states close to the ground states.

We wish now to describe more precise results of a similar character in the context of sums of random independent exponentials, obtained recently by Ben Arous, Bogachev and Molchanov [3], and then state Theorem 3 of this paper which gives more precise information about the region when $L(t) = e^{\gamma c_2 t^{d/(d+2)}}$ with $0 < \gamma < \gamma_2$. In the sequel we will consider infinitely divisible laws described by the Lévy representation of their characteristic function. Let us first set some terminology. We will call a *Lévy spectral function* any function $\mathcal{L}(x):\mathbb{R}/\{0\}\to\mathbb{R}$ which is nondecreasing on $(-\infty,0)$ and on $(0,\infty)$, such that $\lim_{x\to\infty}\mathcal{L}(x)=\lim_{x\to-\infty}\mathcal{L}(x)=0$, and such that for every $\varepsilon>0$ one has $\int_{x:\,\varepsilon\geq|x|>0}x^2\,d\mathcal{L}(x)<\infty$. Now, given a Lévy spectral function $\mathcal{L}$, $\beta\in\mathbb{R}$ and $\sigma>0$, we define $X_{\beta,\sigma,\mathcal{L}}$ as the infinite divisible law whose characteristic function is

$$\phi(u)=\exp\left\{i\beta u-\frac{\sigma^2}{2}u^2+\int_{|x|>0}\left(e^{iux}-1-\frac{iux}{1+x^2}\right)d\mathcal{L}(x)\right\}.$$

In [3] the behavior of sums of the form $S_L(t):=\sum_{k=1}^{L(t)}\exp\{-tX_k\}$ is studied, where $t\in[0,\infty)$, $L(t):[0,\infty)\to\mathbb{R}$ is an increasing-natural-number valued function and $X_1,X_2,\dots$ is a sequence of positive independent identically distributed random variables. In [3] the case where $X_1$ has a Fréchet tail law of parameter $\kappa>0$, so that $P(X_1<1/x)\sim\exp\{-cx^\kappa\}$ for some constant $c>0$, where $P$ is the law of $X_1$, is studied. There it is shown that if $L(t)=\exp\{\gamma H(t)\}$, with $H(t):=\log E(e^{tX_1})$ and $E$ the expectation associated to $P$, then there exist $\gamma_1$ and $\gamma_2$ such that $0<\gamma_1<\gamma_2$ and for which if $\gamma>\gamma_1$, $S_L(t)/E(S_L(t))$ converges to 1 in $P$-probability (law of large numbers) while if $\gamma>\gamma_2$ in addition it is true that $(S_L(t)-E(S_L(t)))/\sqrt{\operatorname{Var}S_L(t)}$, where Var is the variance associated to $P$, converges in distribution to the normal law $\mathcal{N}(0,1)$ (central limit theorem). Moreover, in [3] it is shown that although the central limit theorem fails when $\gamma_1<\gamma\leq\gamma_2$, it is true that $(S_L(t)-E(S_L(t)))/B_1(t)$ converges in law to an infinitely divisible distribution $X_{\beta_1,0,\mathcal{L}_1}$, where $\mathcal{L}_1=c_1/x^a$ for $x>0$ and $\mathcal{L}_1(x)=0$ for $x<0$ with $c_1>0$ and $a=(\frac{\kappa}{\kappa'}\gamma)^{1/\kappa'}$, so that $X_{\beta_1,0,\mathcal{L}_1}$ is a stable law of characteristic exponent $a$. Here, $\kappa'$ is negative conjugate to $\kappa$ (i.e., $\frac{1}{\kappa'}-\frac{1}{\kappa}=1$), while $B_1(t,\gamma):=\exp\{c(\gamma)H(t)\}$, with $H(t):=\log E(e^{tX_1})$ and

$$c(\gamma):=\kappa'\left(\frac{\kappa}{\kappa'}\gamma\right)^{1/\kappa}. \tag{10}$$

Furthermore, they prove that in the region $0<\gamma<\gamma_1$ where the law of large numbers fails, it is true that $S_L(t)/B_1(t)$ converges in law to $X_{\beta_2,0,\mathcal{L}_2}$, where $\beta_2$ is some real number and $\mathcal{L}_2(x)=c_2/x^a$ for $x>0$ and $\mathcal{L}(x)=0$ for $x<0$ with $c_2>0$ and $a=(\frac{\kappa}{\kappa'}\gamma)^{1/\kappa'}$. Thus, the limiting law is stable of characteristic exponent $a$. In [3] a similar analysis is made in the case in which $S_L(t):=\sum_{k=1}^{L(t)}\exp\{tX_k\}$, $X_1,X_2,\dots$ is a sequence of i.i.d. random variables and $P(X>x)\sim\exp\{-cx^\kappa\}$ for some parameter $\kappa>0$ and constant $c>0$.

The results of Theorem 2 might suggest that the averaged survival probability $p^{L(t)}(0,t,w)$ should behave in a way similar to the random exponentials studied in [3], so that properly centered and normalized it converges to stable laws, in the regions of Theorem 1 where the central limit theorem or the law of large numbers fails. The following theorem, however, shows that at least in the one-dimensional case there is convergence to specific infinitely divisible laws which are not stable. In fact, the discrete character of the random walk law manifests itself. Throughout the sequel, given a real number $x$ we will denote by $[x]$ the integer part of $x$ and denote $[x]_- := \lim_{\varepsilon\to 0^-}[x-\varepsilon]$.

THEOREM 3. *Let $d=1$. Let $\alpha=\frac{1}{2}$ and $\alpha'=\frac{1}{3}$ be negative conjugate to $\alpha$ and $L(t) := \exp\{\nu[\frac{\gamma c_2}{\nu}t^{1/3}]_-\}$. The following statements are true.*

(i) Case 3. *If $0<\gamma<\gamma_1$, then we have that*

$$\lim_{t\to\infty}\frac{\sum_{x\in\Lambda_{L(t)}} p(x,t,w)}{s_1(\gamma)t^{1/3}\exp\{-4t\ell_1/[\gamma c_2 t^{1/3}/\nu]_-^2\}} = X_{\beta_1,0,\mathcal{L}}. \tag{11}$$

(ii) Case 4. *If $\frac{2}{3}=\gamma_1<\gamma<\gamma_2=2^{1/3}\frac{2}{3}$, then we have that*

$$\lim_{t\to\infty}\frac{\sum_{x\in\Lambda_{L(t)}}(p(x,t,w)-\langle p(0,t)\rangle)}{s_1(\gamma)t^{1/3}\exp\{-4t\ell_1/[\gamma c_2 t^{1/3}/\nu]_-^2\}} = X_{\beta_2,0,\mathcal{L}}.$$

*In both cases the convergence is in distribution, $s_1(\gamma) := \frac{\gamma c_2}{\ell_1\nu}$, $\beta_1 := \frac{2p}{(1-p)}\times\sum_{k=-\infty}^{\infty}\frac{(1-p)^k}{(1-p)^{k/a_1}+(1-p)^{-k/a_1}}$, $\beta_2 := \frac{2p}{(1-p)}\sum_{k=-\infty}^{\infty}\frac{(1-p)^{k(1+2/a_1)}}{(1-p)^{k/a_1}+(1-p)^{-k/a_1}}$, the Lévy spectral function $\mathcal{L}(x) := -\frac{2p}{1-p}(1-p)^{[-(1/\nu)\log x^{a_1(\gamma)}]}$ for $x>0$ and $\mathcal{L}(x)=0$ for $x<0$ and $a_1 := (\frac{\alpha}{\alpha'}\gamma)^{1/\alpha'}$.*

REMARK 1. Note that the expression for the quantity $a(\gamma)-\gamma$ of Theorem 2 is as the expression for $c(\gamma)$ in the stable limit law cases of [3], in (10). Note that if $a(\gamma) := \alpha'(\frac{\alpha}{\alpha'}\gamma)^{-1/\alpha}+\gamma$, then $(a-\gamma)c_2t^{1/3}\sim\frac{4t\ell_1}{[\gamma c_2 t^{1/3}/\nu]_-^2}$. This is in accordance with the asymptotics described by Theorem 2. On the other hand, the extension of Theorem 3 to higher dimensions requires a better understanding of the tail distribution of the Dirichlet principal eigenvalue of the Laplacian operator on the set $\Lambda_L\cap\mathcal{G}(w)^c$.

REMARK 2. It is possible to prove a one-dimensional theorem analogous to Theorem 3 in the context of Brownian motion on Poissonian obstacles. In this case the limiting laws should be stable.

**3. Moments and correlations.** To prove Theorems 1, 2 and 3, we will need several results describing the asymptotic behavior of some moments and correlations of the field of quenched survival probabilities $\{p(x,t,w) : x\in$ $\mathbb{Z}^d\}$ and truncated versions of them. Here we will prove them, introducing the corresponding notation. This section is divided in three subsections. In the first one we introduce the truncated survival probabilities which will play an important role in reducing some computations to sums of independent random variables. In the next subsection some estimates for the survival probabilities in terms of principal Dirichlet eigenvalues will be proved. Finally, in the third subsection we apply the spectral estimates to compute the logarithmic asymptotics of some quantities depending on moments and correlations of the survival probabilities.

3.1. *Truncated survival probabilities.* Let us first define, given a subset $U \subset \mathbb{Z}^d$, the first exit time $T_U := \inf_{s\geq 0}\{Z_t \notin U\}$ of the random walk $Z.$ from the set $U$. Next, for $x \in \mathbb{Z}^d$, $t \in [0,\infty)$ and $w \in X$, let

$$\tilde{p}_U(x,t,w) := P_x(\tau(w) > t, T_U \geq t). \tag{12}$$

Given $x \in \mathbb{Z}^d$ and $t \geq 0$, in the particular case in which $U = \Lambda(x, at)$, we will use the notation $\tilde{p}_a(x,t,w)$ instead of $\tilde{p}_U(x,t,w)$ and $T_{at}$ instead of $T_U$. We will call the probabilities $\tilde{p}_a(x,t,w)$ the *truncated quenched survival probabilities* at scale $a$ at time $t$ for a random walk starting from $x$ and refer to $\langle \tilde{p}_a(x,t)\rangle$ as the *truncated annealed survival probabilities* at scale $a$ at time $t$ for a random walk starting from $x$. Furthermore, we will denote the sets $\{\tilde{p}_a(x,t,w) : x \in \mathbb{Z}^d\}$, $\{\langle \tilde{p}_a(x,t)\rangle : x \in \mathbb{Z}^d\}$ as the *field of truncated quenched survival probabilities* and the *field of truncated annealed survival probabilities* at scale $a$, respectively. These quantities will be useful to approximate the field of quenched survival probabilities by an independent field defined on boxes. This will eventually be applied for the proof of the law of large numbers and central limit theorems stated in parts (ii) and (iii) of Theorem 1. We will prove results describing the long-time behavior of the moments $\langle p(0,t)^r\rangle$ and $\langle \tilde{p}_a(0,t)^r\rangle$, for $r$ real and $r \geq 1$, and of the correlations $\langle p(x,t), p(y,t)\rangle := \langle p(x,t)p(y,t)\rangle - \langle p(x,t)\rangle\langle p(y,t)\rangle$ and $\langle \tilde{p}_a(x,t), \tilde{p}_a(y,t)\rangle := \langle \tilde{p}_a(x,t)\tilde{p}_a(y,t)\rangle - \langle \tilde{p}_a(x,t)\rangle\langle \tilde{p}_a(y,t)\rangle$, where the distance $|x-y|$, between $x$ and $y$, might possibly depend on time. Many of the computations will be done for the survival probabilities regarded as a field, in contrast to the usual estimates that can be found in the literature. We do not consider our estimations novel, however, since they still are very much in the spirit of standard ones. We begin with the following elementary lemma.

LEMMA 1. *Let $a > 0$ and consider the field of quenched survival probabilities $\{p(x,t,w) : x \in \mathbb{Z}^d\}$ and of truncated quenched survival probabilities at scale $a > 0$. Then for every $x \in \mathbb{Z}^d$ and $t \geq 0$,*

$$|p(x,t,w) - \tilde{p}_a(x,t,w)| \leq k_1(a,d)e^{-k_2(a,d)t}, \tag{13}$$

*where* $k_1(a,d) := e^{-\sinh^{-1}(ad)}$ *and* $k_2(a,d) := \frac{1}{d}J(ad)$ *are positive constants with* $J(x) := x\sinh^{-1} x - \sqrt{1+x^2}+1$ *defined for* $x \geq 0$.

PROOF. Note that

$$p(x,t,w) \leq \tilde{p}_a(x,t,w) + P_x(T_{at} < t).$$

Note that the event $T_{at} < t$ happens if and only if one of the $d$ coordinates of the random walk increases by $[at]$ in a time smaller than $t$. Hence,

$$P_x(T_{at} < t) \leq dP\Big(\sup_{0\leq s\leq t} X_s \geq [at]\Big),$$

where $X.$ is a one-dimensional random walk of total jump rate $\frac{1}{d}$. Now, by the reflection principle, the right-hand side of the above inequality can be upper-bounded by $2P(X_t \geq [at])$, where $P$ is the corresponding law. Now $X_t = N_t - M_t$, where $N_t$ and $M_t$ are independent Poisson processes each one of jump rate $1/(2d)$. Therefore, since $[at] \geq at-1$, we have

$$P_x(T_{at} < t) \leq 2dP(N_t - M_t \geq at - 1), \tag{14}$$

where we use for the sake of clarity again $P$ for the joint law of the two Poisson processes. Now, if we call $E$ the expectation with respect to one of these Poisson processes, say $N_t$, we see by the Chebyshev inequality that the probability in the right-hand side of (14) is bounded by $E(e^{N_t\mu}) \cdot E(e^{-N_t\mu})/e^{(at-1)\mu} = e^{(\cosh\mu - 1)t/2d-(at-1)\mu}$, where $\mu > 0$ is arbitrary. Let us remark that this upper bound is optimized for a positive value of $\mu$. Substituting it we obtain

$$P_x(T_{at} < t) \leq k_1(a,d)e^{-k_2(a,d)t},$$

where $k_1(a,d) := e^{-\sinh^{-1}(ad)}$ and $k_2(a,d) := \frac{1}{d}J(ad)$ with $J(x) : [0,\infty) \to [0,\infty)$ defined as $J(x) := x\sinh^{-1} x - \sqrt{1+x^2}+1$. Note that $J'(x) = \log(x + \sqrt{1+x^2})$ is positive for $x > 0$. Hence $k_1(a,d)$ and $k_2(a,d)$ are strictly positive for $a > 0$. $\square$

Estimate (13) tells us how much we lose when replacing the quenched field by the truncated one. It will turn out that this exponential error will be negligible for most of our purposes.

3.2. *Spectral estimates.* We will now proceed to compute bounds for the survival probabilities in terms of principal Dirichlet eigenvalues of a discrete Laplacian operator. Given a subset $U \subset \mathbb{Z}^d$, we define the normalized discrete Laplacian operator by its action on functions $f : \mathbb{Z}^d \to \mathbb{R}$ which vanish outside $U$ [i.e., $f(x) = 0$ for $x \notin U$] as

$$\Delta_d f(x) := \frac{1}{2d}\sum_{e\in\mathcal{B}}(f(x+e) - f(x)), \tag{15}$$

where $\mathcal{B}$ is the union of the elements of the basis of the free abelian group $\mathbb{Z}^d$ and its inverses. Note that if we define $L^2(U) := \{f : \sum_{x\in\mathbb{Z}^d} f^2(x) < \infty f(x) = 0 \text{ if } x \notin U\}$, $\Delta_d$ is a self-adjoint operator on the Hilbert space $L^2(U)$ endowed with the inner product $(f,g) := \sum_{x\in\mathbb{Z}^d} f(x)g(x)$ for $f, g \in L^2(U)$. We can then consider the set $\{\lambda_n(U) : n \in \mathcal{U}\}$ of eigenvalues of $\Delta_d$ in $L^2(U)$ in increasing order, where $\mathcal{U}$ is the index set. Remark that the cardinality of this index set is finite for finite $U$ and at most countable. We will denote by $\{\psi_n^U : n \in \mathcal{U}\}$ the corresponding normalized eigenfunctions [i.e., $\sum_{x\in\mathbf{Z}^d} \psi_n^U(x)^2 = 1$]. Let $r \geq 0$. In the particular case in which $U = \Lambda(x,r) \cap \mathcal{G}^c(w) =: \Lambda(x,r,w)$ we will employ the notation $\{\lambda_n(x,r,w)\}$ instead of $\{\lambda_n(U) : n \in \mathcal{U}\}$ and $\{\psi_n^{x,r,w}\}$ instead of $\{\psi_n^U : n \in \mathcal{U}\}$. Furthermore, in the sequel, $\mathbf{1}_U$ will denote the indicator function of the set $U$. We first begin with the following upper bound.

LEMMA 2. *Consider the field of quenched survival probabilities $\{p(x,t,w) : x \in \mathbb{Z}^d\}$. Then the following statements are true.*

(i) *For every real $a > 0$ it is true that*

$$p(x,t,w) \leq k_1(a,d)e^{-k_2(a,d)t} + (2at+1)^{d/2} e^{-\lambda_0(x,at,w)t}, \tag{16}$$

*where $k_1(a,d) := 2d \sinh^{-1}(ad)$ and $k_2(a,d) := \frac{1}{d}J(ad)$ are positive constants with $J(x)$ defined as in Lemma 1 by $J(x) := x\sinh^{-1} x - \sqrt{1+x^2} + 1$ for $x \geq 0$.*

(ii) *For every finite subset $U \subset \mathbb{Z}^d$,*

$$\frac{1}{|U|}\sum_{z\in U} p(z,t,w) \geq \frac{1}{|U|} e^{-\lambda_0(U\cap\mathcal{G}^c(w))t}. \tag{17}$$

PROOF OF LEMMA 2(i). Note that by Lemma 1 we have that $p(x,t,w) \leq k_1(a,d)e^{-k_2(a,d)t} + \tilde{p}_a(x,t,w)$. Therefore it is enough to estimate the truncated probability at scale $a$, $\tilde{p}_a(x,t,w)$. First remark the following expansion in terms of the eigenvalues $\{\lambda_n(x,at,w)\}$ and the corresponding eigenfunctions $\{\psi_n^{x,at,w}\}$:

$$\tilde{p}_a(x,t,w) = \sum_{n\in\mathcal{U}} e^{-t\lambda_n(x,at,w)} \psi_n^{x,at,w}(x)(\psi_n^{x,at,w}, \mathbf{1}_A), \tag{18}$$

where $A := \Lambda(x,at,w)$. Now, by the Cauchy–Schwarz inequality we see that the right-hand side of (18) is upper-bounded by $e^{-t\lambda_0(x,at,w)}(\sum_{n\in\mathcal{U}}(\psi_n^{x,at,w}, \mathbf{1}_{x,w})^2 \sum_{n\in\mathcal{U}}(\psi_n^{x,at,w}, \mathbf{1}_A)^2)^{1/2}$ which in turn is upper-bounded by $e^{-t\lambda_0(x,at,w)} \times \sqrt{|A|}$, where $\mathbf{1}_{x,w}(y)$ equals 1 if $y = x \in \mathcal{G}^c(w)$ and 0 otherwise. Using the fact that $|A| = |\Lambda(x,at,w)| \leq ([2at]+1)^d$, we conclude the proof. $\square$

Proof of Lemma 2(ii). First note the trivial inequality $p(z,t,w) \geq \tilde{p}_U(z,t,w)$, where $\tilde{p}_U(z,t,w)$ is defined in (12). Also remark the following eigenfunction expansion: $\tilde{p}_U(z,t,w) = \sum_{n\in\mathcal{U}} e^{-\lambda_n(U\cap\mathcal{G}^c(w))t}\psi_n^{U\cap\mathcal{G}^c(w)}(z) \times (\psi_n^{U\cap\mathcal{G}^c(w)}, \mathbf{1}_U)$. Therefore we can see that

$$\begin{aligned}\frac{1}{|U|}\sum_{z\in U} p(z,t) &\geq \frac{1}{|U|}e^{-\lambda_0(U\cap\mathcal{G}^c(w))t}(\psi_0^{U\cap\mathcal{G}^c(w)}, \mathbf{1}_U)^2 \\ &\geq \frac{1}{|U|}e^{-\lambda_0(U\cap\mathcal{G}^c(w))t}\sum_{z\in U}(\psi_0^{U\cap\mathcal{G}^c(w)})^2(z) \\ &= \frac{1}{|U|}e^{-\lambda_0(U\cap\mathcal{G}^c(w))t},\end{aligned}$$

where we have used in the second to last inequality the fact that $\psi_0^{U\cap\mathcal{G}^c(w)}(x) \geq 0$ and in the last inequality the normalization condition $\sum_{z\in U}(\psi_0^{U\cap\mathcal{G}^c(w)})^2(z) = 1$.

□

3.3. *Asymptotics of moments and correlations.* Let us now apply the previous results to estimate quantities such as some moments and correlations of the fields of quenched and truncated quenched survival probabilities.

We begin writing a couple of elementary inequalities that will be repeatedly used throughout the paper. Let $n$ be an arbitrary natural number and let $a_1,\dots,a_n$ be arbitrary real numbers. If $r \geq 1$ we have by Jensen's inequality

$$\left|\sum_{i=1}^n a_i\right|^r \leq n^{r-1}\sum_{i=1}^n |a_i|^r, \tag{19}$$

while if $0 \leq r \leq 1$ we have

$$\left|\sum_{i=1}^n a_i\right|^r \leq \sum_{i=1}^n |a_i|^r. \tag{20}$$

Our first result is a pair of lemmas for the field of quenched and truncated quenched survival probabilities, respectively.

Lemma 3. *Consider the field of quenched survival probabilities,* $\{p(x,t,w): x\in\mathbb{Z}^d\}$. *Let* $x,y\in\mathbb{Z}^d$. *Then:*

(i) *For every real* $r > 0$ *it is true that*

$$\log\langle p(0,t)^r\rangle \sim -c_2(d,p)(rt)^{d/(d+2)}. \tag{21}$$

(ii) *For every real $r>0$ it is true that*

$$\log\langle|p(x,t)-\langle p(x,t)\rangle|^r\rangle \sim -c_3(d,p,r)t^{d/(d+2)},$$

*where* $c_3(d,p,r):=c_2(d,p)\min\{r,r^{d/(d+2)}\}$.

(iii) *For every $t\geq 0$ it is true that*

$$\langle p(x,t),p(y,t)\rangle \geq 0.$$

(iv) *Let $\{U_t:t>0\}$ be a collection of subsets of the lattice $\mathbb{Z}^d$ indexed by $t>0$. Assume that there is an $a>0$ such that $|U_t|\sim|U_{a,t}|$ as $t\to\infty$, where $U_{a,t}:=\{x\in U_t:\operatorname{dist}(x,U_t^c)\geq 2at\}$. Then,*

$$\operatorname{Var}_\mu\left(\sum_{x\in U_t}p(x,t)\right)\sim|U_t|\left(\sum_{x\in\mathbb{Z}^d}\langle p(0,t),p(x,t)\rangle\right). \tag{22}$$

*Furthermore,*

$$\log\left(\sum_{x\in\mathbb{Z}^d}\langle p(0,t),p(x,t)\rangle\right)\sim -c_2(d,p)(2t)^{d/(d+2)}. \tag{23}$$

REMARK 3. Part (i) of Lemma 3 in the case of integer-valued $r\geq 1$ and for Brownian motion is stated in Section 4.5 of [12].

LEMMA 4. *Consider the field of truncated quenched survival probabilities at scale $a>0$, $\{\tilde p_a(x,t,w):x\in\mathbb{Z}^d\}$. Let $x,y\in\mathbb{Z}^d$. Then:*

(i) *For every real $r>0$ it is true that*

$$\log\langle\tilde p_a(0,t)^r\rangle\sim -c_2(d,p)(rt)^{d/(d+2)}.$$

(ii) *For every real $r>0$ it is true that*

$$\log\langle|\tilde p_a(x,t)-\langle\tilde p_a(x,t)\rangle|^r\rangle\sim -c_3(d,p,r)(rt)^{d/(d+2)},$$

*where* $c_3(d,p,r):=c_2(d,p)\min\{r,r^{d/(d+2)}\}$.

(iii) *For every $t\geq 0$ it is true that*

$$\langle\tilde p_a(x,t),\tilde p_a(y,t)\rangle + k_1(a,d)e^{-k_2(a,d)t}\geq 0, \tag{24}$$

*where $k_1(a,d):=2d\sinh^{-1}(ad)$ and $k_2(a,d):=\frac{1}{d}J(ad)$ are positive constants with $J(x)$ defined as in Lemma 1 by $J(x):=x\sinh^{-1}x-\sqrt{1+x^2}+1$ for $x\geq 0$.*

(iv) *Let $\{U_t:t>0\}$ be a collection of subsets of the lattice $\mathbb{Z}^d$ indexed by $t>0$. Assume that $|U_t|\sim|U_{a,t}|$ as $t\to\infty$, where $U_{a,t}:=\{x\in U_t:\operatorname{dist}(x,U_t^c)\geq 2at\}$. Then,*

$$\operatorname{Var}_\mu\left(\sum_{x\in U_t}\tilde p_a(x,t)\right)\sim|U_t|\left(\sum_{x\in\mathbb{Z}^d}\langle\tilde p_a(0,t),\tilde p_a(x,t)\rangle\right). \tag{25}$$

*Furthermore,*

$$\log\left(\sum_{x\in\mathbb{Z}^d}\langle\tilde{p}(0,t),\tilde{p}_a(x,t)\rangle\right)\sim -c_2(d,p)(2t)^{d/(d+2)}. \tag{26}$$

The proofs of parts (i), (ii) and (iii) of Lemma 4 are analogous to the proofs of the corresponding parts of Lemma 3. Therefore, we will prove parts (i), (ii) and (iii) of Lemma 3 and subsequently, part (iv) of Lemma 4 and Lemma 3.

PROOF OF LEMMA 3(i). Let $x\in\mathbb{Z}^d$ and $t\geq 0$. Note that by part (ii) of Lemma 2 with $U=\Lambda(x,t)$ in (17), it is true that

$$\frac{1}{|\Lambda(x,t)|}\sum_{z\in\Lambda(x,t)}p(z,t)\geq\frac{1}{(2t+1)^d}e^{-\lambda_0(x,t,w)t}. \tag{27}$$

Taking expectations on both sides of (27) we conclude that

$$\langle p(0,t)\rangle\geq\frac{1}{(2t+1)^d}\langle e^{-\lambda_0(0,t,w)t}\rangle, \tag{28}$$

where we have used the translation invariance of the measure $\mu$. On the other hand, choosing $a=r$ in (16), we see that

$$p(0,t)\leq k_1(r,d)e^{-k_2(r,d)t}+(2rt+1)^{d/2}e^{-\lambda_0(0,rt,w)t}. \tag{29}$$

From Jensen's inequality (19) with $n=2$ and (19) applied to (29) we deduce that for every $r\geq 0$,

$$\langle p(0,t)^r\rangle\leq 2^r(k_1^re^{-k_2rt}+(2rt+1)^{rd/2}\langle e^{-\lambda_0(0,rt,w)rt}\rangle). \tag{30}$$

Combining (28), (29) and (30), with the well-known fact that

$$\lim_{t\to\infty}\frac{\log\langle p(0,t)\rangle}{c_2(d,p)t^{d/(d+2)}}=-1$$

(see [2] or [5]), we deduce that for every $r>0$,

$$\limsup_{t\to\infty}\frac{\log\langle p(0,t)^r\rangle}{c_2(d,p)(rt)^{d/(d+2)}}\leq -1. \tag{31}$$

The lower bound,

$$\liminf_{t\to\infty}\frac{\log\langle p(0,t)^r\rangle}{c_2(d,p)(rt)^{d/(d+2)}}\geq -1, \tag{32}$$

can be proved following the standard approach of [4], [5] or [2] (see also [11]). In what follows we recall such procedure in our context. First we note that for every real $R\geq 0$ we have

$$\langle p(0,t)^r\rangle\geq\left\langle A_RP_0\left(\max_{0\leq s\leq t}|Z_s|\leq R\right)^r\right\rangle,$$

where $A_R$ is the indicator function of the event that the ball $B(0,R)$ does not contain any site which is an obstacle. Now, an estimation of the random walk probability from below using an eigenvalue expansion and retaining only the term with principal Dirichlet eigenvalue, and an optimization with respect to $R$ enable us to deduce (32). Inequalities (31) and (32) together imply the statement (21). □

Proof of Lemma 3(ii). By (19) and (20) note that $|p(0,t)-\langle p(0,t)\rangle|^r \le 2^r(p(0,t)^r + \langle p(0,t)\rangle^r)$. This and part (i) of Lemma 3 then imply the following upper bound for $r>0$:

$$\limsup_{t\to\infty} \frac{\log\langle |p(0,t)-\langle p(0,t)\rangle|^r\rangle}{c_3(d,p,r)t^{d/(d+2)}} \le -1.$$

On the other hand,

$$\begin{aligned}\langle |p-\langle p\rangle|^r\rangle &\ge \langle |p-\langle p\rangle|^r \mathbf{1}_{p>2\langle p\rangle}\rangle \ge \left\langle p^r \left|1-\frac{\langle p\rangle}{p}\right|^r \mathbf{1}_{p>2\langle p\rangle}\right\rangle \\ &\ge \frac{1}{2^r}\langle p^r \mathbf{1}_{p>2\langle p\rangle}\rangle \ge \frac{1}{2^r}(\langle p^r\rangle - (2\langle p\rangle)^r).\end{aligned}$$

From here we conclude that whenever $r>1$ it is true that

$$\liminf_{t\to\infty} \frac{\log\langle |p(0,t)-\langle p(0,t)\rangle|^r\rangle}{c_3(d,p,r)t^{d/(d+2)}} \ge -1. \tag{33}$$

On the other hand, note that if $a>0$ and $b>0$, then whenever $0\le r\le 1$, it is true that $|a-b|^r \ge |a^r - b^r|$. Thus, if $0\le r\le 1$ we have that $|p-\langle p\rangle|^r \ge \langle p\rangle^r - p^r$. This clearly implies that the lower bound (33) is also valid for $0\le r<1$. If $r=1$, note from the inequalities $|p-\langle p\rangle|^{1+\varepsilon} \le |p-\langle p\rangle| \le |p-\langle p\rangle|^{1-\varepsilon}$ valid for $\varepsilon>0$ that $-(1+\varepsilon)^{d/(d+2)} \le \liminf_{t\to\infty} \frac{\log\langle |p(0,t)-\langle p(0,t)\rangle|\rangle}{c_2(d,p)t^{d/(d+2)}} \le -(1-\varepsilon)$. Letting $\varepsilon\to 0$ we conclude that the lower bound (33) is also valid for $r=1$. □

Proof of Lemma 3(iii). The case $x=y$ is trivial. So let us assume that $x\ne y$. Denote by $E_x$ the expectation of the continuous-time simple random walk $Z.$ of total jump rate 1 starting from site $x\in\mathbf{Z}^d$ and of law $P_x$. Define the Wiener sausage $W_x(t)$ at time $t$ as the sites visited between time 0 and $t$ for the random walk $Z.$ starting from $x$ so that $W_x(t) := \{z\in\mathbf{Z}^d : z = Z_s \quad \text{for some time } s\in[0,t]\}$. Now note that $\langle p(x,t)\rangle = E_x(e^{-\nu|W_x(t)|})$ and $\langle p(x,t)p(y,t)\rangle = E_{x,y}(e^{-\nu|W_x(t)\cup W_y(t)|})$, where $E_{x,y} := E_x\otimes E_y$ denotes expectation with respect to independent random walks with laws $P_x$ and $P_y$ and $\nu := |\log(1-p)|$. Hence,

$$\langle p(x,t),p(y,t)\rangle = E_{x,y}(e^{-\nu(|W_x(t)+W_y(t)|)}(e^{\nu|W_x(t)\cap W_y(t)|}-1)). \tag{34}$$

Part (ii) follows immediately. $\square$

Proofs of Lemma 3(iv) and Lemma 4(iv). Note that due to the fact that the truncated survival probabilities $\tilde{p}_a(x,t)$ and $\tilde{p}_a(y,t)$ are independent for $\|x-y\|>2at$, we have that

$$\operatorname{Var}_\mu\left(\sum_{x\in U_t}\tilde{p}_a(x,t)\right)=\sum_{x,y\in U_t}\langle\tilde{p}_a(x,t),\tilde{p}_a(y,t)\rangle$$
$$=\sum_{x,y\in U_t\colon\|x-y\|\le 2at}\langle\tilde{p}_a(x,t),\tilde{p}_a(y,t)\rangle.$$

Now, the rightmost member of the above equalities is bounded above by $|U_t|\sum_{y\colon\|y\|\le 2at}\langle\tilde{p}_a(0,t),\tilde{p}_a(y,t)\rangle$ which in turn gives the upper bound

$$\operatorname{Var}_\mu\left(\sum_{x\in U_t}\tilde{p}_a(x,t)\right)\le|U_t|\sum_{y\in\mathbf{Z}^d}\langle\tilde{p}_a(0,t),\tilde{p}_a(y,t)\rangle. \tag{35}$$

On the other hand, in a similar way we can show that

$$\operatorname{Var}_\mu\left(\sum_{x\in U_t}\tilde{p}_a(x,t)\right)\ge|U_{a,t}|\sum_{y\in\mathbf{Z}^d}\langle\tilde{p}_a(0,t),\tilde{p}_a(y,t)\rangle. \tag{36}$$

But, since $|U_t|\sim|U_{a,t}|$, inequalities (35) and (36) prove (25). To prove (26), note that

$$\sum_{x\in\mathbf{Z}^d}\langle\tilde{p}_a(0,t),\tilde{p}_a(x,t)\rangle\ge\operatorname{Var}_\mu(\tilde{p}_a(0,t)). \tag{37}$$

However, it is true that $\sum_x\langle\tilde{p}_a(0,t),\tilde{p}_a(x,t)\rangle=\sum_{x\colon\|x\|\le 2at}\langle\tilde{p}_a(0,t),\tilde{p}_a(x,t)\rangle$. This combined with the inequality $\langle\tilde{p}_a(0,t),\tilde{p}_a(x,t)\rangle\le\operatorname{Var}_\mu(\tilde{p}_a(0,t))$ for $x\in\mathbf{Z}^d$, shows that

$$\sum_{x\in\mathbf{Z}^d}\langle\tilde{p}_a(0,t),\tilde{p}_a(x,t)\rangle\le(3at)^d\operatorname{Var}_\mu(\tilde{p}_a(0,t)). \tag{38}$$

Inequalities (37) and (38) together with the asymptotic behavior (24) with $r=2$, imply (26). We continue now with the proof of part (iv) of Lemma 3. Let $a>0$ be such that $|U_t|\sim|U_{a,t}|$. Note that from (34) we have

$$\langle\tilde{p}_a(x,t),\tilde{p}_a(y,t)\rangle\le\langle p(x,t),p(y,t)\rangle, \tag{39}$$

for every $x,y\in\mathbb{Z}^d$ and $t\ge 0$. This, and a calculation similar to the one which leads to (35), enables us to conclude that

$$\operatorname{Var}_\mu\left(\sum_{x\in|U_t|}\tilde{p}_a(x,t)\right)\le\operatorname{Var}_\mu\left(\sum_{x\in|U_t|}p(x,t)\right)\le|U_t|\sum_{y\in\mathbf{Z}^d}\langle p(0,t),p(y,t)\rangle.$$

Hence, to complete the proof of (22) and (23) it is enough to show that

$$\sum_{y\in\mathbf{Z}^d}\langle p(0,t),p(y,t)\rangle \sim \sum_{y\in\mathbf{Z}^d}\langle \tilde p_a(0,t),\tilde p_a(y,t)\rangle, \tag{40}$$

applying part (iv) of Lemma 4. Now,

$$\begin{aligned}\sum_{y\in\mathbb{Z}^d}\langle p(0,t),p(y,t)\rangle &= \sum_{y\,:\,\|y\|\le 2at}\langle p(0,t),p(y,t)\rangle\\ &\quad+\sum_{y\,:\,\|y\|>2at}\langle p(0,t),p(y,t)\rangle.\end{aligned}\tag{41}$$

Via Lemma 1, the last term can be shown to be bounded by

$$\sum_{y\,:\,\|y\|>2at} 2k_1\left(\frac{\|y\|}{2t},d\right)e^{-k_2(\|y\|/2t,d)t},$$

which is exponentially small in $t$. A second application of the same result lets us conclude that $\sum_{y\,:\,\|y\|\le 2at}\langle p(0,t),p(y,t)\rangle \sim \sum_{y\in\mathbb{Z}^d}\langle \tilde p_a(0,t),\tilde p_a(y,t)\rangle$. These remarks, together with the monotonicity (39) and equality (41), prove (40). □

**4. The quenched and intermediate asymptotics.** In this section we will prove parts (i) and (i′) of Theorem 1. We will indirectly make use of the method of enlargement of obstacle through some standard estimates on the almost sure asymptotic behavior of principal Dirichlet eigenvalues and on the geometry of the obstacle set $\mathcal{G}(w)$.

4.1. *Quenched asymptotics.* Let us begin proving part (i). First we show that if $p<1-p_c$, then $\mu$-a.s.,

$$\liminf_{t\to\infty}\frac{(\log t)^{2/d}}{c_1 t}\log p^{L(t)}(0,t,w)\ge -1.$$

Now, this is a simple consequence of the inequality $p^{L(t)}(0,t,w)\ge \frac{1}{(2t+1)^d}\,p(0,t,w)$, where we used that $L(t)\le t$, and the fact that $\mu$-a.s. on the event $C(w)$ that the origin belongs to an infinite trap free cluster it is true that $\liminf_{t\to\infty}\frac{(\log t)^{2/d}}{c_1 t}\log p(0,t,w)\ge -1$ (see [2] or Theorem 4.5.1 on page 196 of [12] for a Brownian motion version of this result).

Let us now prove that $\mu$-a.s.,

$$\limsup_{t\to\infty}\frac{(\log t)^{2/d}}{c_1 t}\log p^{L(t)}(0,t,w)\le -1. \tag{42}$$

Now note that for every $x\in\Lambda_{L(t)}$ it is true that $\lambda_0(x,t,w)\ge\lambda_0(0,L(t)+t,w)\ge\lambda_0(0,2t,w)$, where in the last inequality we have used the fact that $L(t) \le t$. Hence, for every $x \in \Lambda_{L(t)}$ we have that $e^{-\lambda_0(x,t,w)t} \le e^{-\lambda_0(0,2t,w)t}$ so that by (16) with $a = 1$, we see that

$$p^{L(t)}(0,t,w) \le k_1(1,d)e^{-k_2(1,d)t} + (2t+1)^{d/2}e^{-\lambda_0(0,2t,w)t}. \tag{43}$$

Finally, standard enlargement of obstacle asymptotic estimates for $\lambda_0(0,2t,w)$ (see [1]) tell us that $\mu$-a.s. it is true that $\lim_{t\to\infty}(\log t)^{d/2}\lambda_0(0,2t,w) = c_1(d,p)$. This fact, combined with the upper bound (43), implies (42).

4.2. *Intermediate asymptotics.* We now proceed to prove part (i′) of Theorem 1. We first show that

$$\liminf_{t\to\infty} \frac{(\log L(t))^{2/d}}{c_1 t} p^{L(t)}(0,t,w) \ge -1. \tag{44}$$

Let us recall that $R_0 := (\frac{d}{w_d \nu})^{1/d}$. The principal ingredient in the proof of the limit (44) will be the following statement which is a random walk version of the Brownian motion context Lemma 4.5.2 of [12]. If $L(t) \ge t$, $\mu$-a.s. eventually in $t$, the following event occurs:

$$\bigcup_{x\in\Lambda_{L'(t)}} \{B(x, R(\log L(t))^{1/d}) \cap (\mathcal{G}(w))^c = \varnothing\}, \tag{45}$$

where $R := R_0 - \varepsilon(t)$ with $\varepsilon(t) = 1/(\log L(t))^{1/d}$ for $d \ge 2$ while $\varepsilon(t) := \frac{4R_0}{3}(\log\log L(t))/\log L(t)$ if $d = 1$, and $L'(t) := L(t) - R_0(\log L(t))^{1/d}$. In other words, $\mu$-a.s. eventually in $t$, there exists a Euclidean ball of radius $R(\log L(t))^{1/d}$ contained in $\Lambda_{L(t)}$ which has no obstacles. To prove this, first note that the box $\Lambda_{L(t)}$ contains at least $[(2L+1)/(2R_0(\log L)^{1/d})]^d$ disjoint boxes of side $2R_0(\log L(t))^{1/d}$. Now, the probability of the event (45) is smaller than the probability that some of these boxes contain a Euclidean ball of radius $R(\log L(t))^{1/d}$ which has no obstacles. This quantity is smaller than $(1 - \exp\{-\nu w_d R^d \log L\})^{(2L+1)^d/((2R_0)^d \log L)-1}$, which in turn is bounded by $e \cdot \exp\{-e^{-\nu w_d R^d \log L}\frac{(2L+1)^d}{(2R_0)^d \log L}\}$. Using the inequality $(R_0 - \varepsilon(t))^d \ge R_0^d - dR_0^{d-1}\varepsilon(t)$, we can bound this quantity by $\exp\{-\frac{3^d}{(2R_0)^d}(\log L)^{-1}e^{c(\varepsilon)}\}$, where $c(\varepsilon) := d\nu\, w_d R_0^{d-1}\varepsilon^{-(d-1)}$ when $d \ge 2$, while $c(\varepsilon) := \frac{4}{3}\log\log L$ when $d = 1$. For $L$ large enough this is smaller than $L^{-d}$. An application of the Borel–Cantelli lemma proves our claim (45). On the other hand, by (17) with $U = \Lambda_{L(t)}$, we have that

$$p^{L(t)}(0,t,w) \ge \frac{1}{(2L+1)^d} e^{-\lambda_0(0,L,w)t}. \tag{46}$$

Now, by (45), the monotonicity of the principal Dirichlet eigenvalue with respect to the partial order of inclusion of sets, and translation invariance, it is true that $\lambda_0(0,L,w) \le \lambda_0(B(0,R(\log L)^{1/d}))$. Hence, from (46) we conclude that $p^{L(t)}(0,t,w) \ge \frac{1}{(2L+1)^d} e^{-\lambda(B(0,R(\log L)^{1/d}))t}$. Now, using the fact that $\lim_{L\to\infty} \lambda(B(0,R(\log L)^{1/d}))(\log L)^{2/d} = \ell_d/R_0^2$, where we recall that $\ell_d$ denotes the principal Dirichlet eigenvalue of the continuous Laplacian operator on the ball of unit radius (see [2]), and using that $\frac{t}{(\log L)^{2/d}} \gg \log L$ (which is a consequence of the hypothesis $\log L \ll t^{d/(d+2)}$) we conclude the lower bound (44).

Let us now show that $\mu$-a.s.,

$$\limsup_{t\to\infty} \frac{(\log L(t))^{2/d}}{c_1 t} \log p^{L(t)}(0,t,w) \le -1. \tag{47}$$

The proof of (47) is very similar to the proof of the upper bound (42). This time note that from (16) with $a=1$ we deduce that

$$p^{L(t)}(0,t,w) \le k_1(1,d)e^{-k_2(1,d)t} + (2t+1)^{d/2} e^{-\lambda_0(0,t+L,w)t}. \tag{48}$$

As in part (i) the first term is negligible. On the other hand, since for $L(t) \gg t$, we have that $\mu$-a.s. $\lambda_0(\Lambda_{L(t)+t}(w)) \sim \lambda_0(\Lambda_{L(t)}(w))$, the desired upper bound (47) follows.

**5. The critical, annealed and Gaussian asymptotics.** In this section we will prove parts (ii) and (iii) of Theorem 1, concerning the law of large numbers and central limit theorem, and Theorem 2. The proofs are based on a Dirichlet–Neumann type partition analysis. This enables us to arrive up to logarithmically sharp lower bounds for the growth of $L(t)$ so that the law of large numbers and the central limit theorem in Theorem 1 are valid, and to determine the exact rate of growth of the quantities $p^L$ and $p^L - \langle p\rangle$ in the corresponding complementary regions. Our method is based on making a partition of the box $\Lambda_L$ in a collection of smaller boxes indexed by a set $\mathcal{I}$, $\{\Lambda_{\mathbf{i}} : \mathbf{i} \in \mathcal{I}\}$. In some computations, this collection will in turn be subdivided in $2^d$ disjoint collections which essentially decouple the dependence of the random variables appearing in the sums over $\Lambda_L$. In others, it will be necessary to define a strip set whose contribution will be negligible, and such that the set of partition boxes intersected with the complement of the strip set becomes an independent set. This section is divided in six subsections. First we define the decomposition of the box $\Lambda_L$ in the described partition boxes and strip set. In the second subsection we collect several technical results that will be used repeatedly. In the third subsection we prove the law of large numbers stated in (1) up to scales $L(t) \ge \exp\{\gamma \frac{c_2}{d} t^{d/(d+2)}\}$ with $\gamma > \gamma_1$. In the fourth subsection we prove the central limit theorem stated in (4), up to scales $L(t) \ge \exp\{\gamma \frac{c_2}{d} t^{d/(d+2)}\}$ with $\gamma > \gamma_2$. In the fifth subsection we prove part (i) of Theorem 2 obtaining as a corollary (3) stating that the value $\gamma_1$ is the smallest possible $\gamma$ for the validity of the law of large numbers (1). Similarly in the last subsection we prove part (ii) of Theorem 2 obtaining as a corollary (4) stating that the value $\gamma_2$ is the smallest possible $\gamma$ for the validity of the central limit theorem (4).

5.1. *Partition analysis.* Let $L$ be some natural number and consider the corresponding box $\Lambda_L = \{x \in \mathbb{Z}^d : \|x\| \le L\}$. Here we will define two related but different kinds of partitions of $\Lambda_L$. The first one shows that $\Lambda_L$ can be decomposed in disjoint *partition boxes* $\{\Lambda'_{\mathbf{i}} : \mathbf{i} \in \mathcal{I}\}$, indexed by some set $\mathcal{I}$, so that $\Lambda_L = \bigcup_{\mathbf{i}} \Lambda'_{\mathbf{i}}$. The second one defines a partition of $\Lambda_L$ in a *strip set* and *main boxes* $\{\Lambda''_{\mathbf{i}} : \mathbf{i} \in \mathcal{I}\}$. In the first case, the index set $\mathcal{I}$ will be partitioned in disjoint subsets $\{\mathcal{I}_K : K \in \mathcal{K}\}$, where the cardinality of $\mathcal{K}$ is $2^d$, in such a way that for each $K \in \mathcal{K}$ any pair of elements of the collection of partition boxes $\{\Lambda'_{\mathbf{i}} : \mathbf{i} \in \mathcal{I}_K\}$ is at a large Euclidean distance. This property will enable us to approximately decouple sums of the form $\sum_{x\in\Lambda_L} p$ or $\sum_{x\in\Lambda_L}(p - \langle p\rangle)$ in a finite and constant number of sums of independent random variables. In the second partition case, it turns out that the survival probabilities corresponding to sites in the strip set have a total sum which is negligible, while the main boxes happen to be essentially independent. To proceed we will need to introduce some notation defining the corresponding scales and subsets.

Our first parameter is a natural number $L'$ smaller than or equal to $L$. Throughout the sequel $L'$ will be called the *mesoscopic scale*. By the division algorithm, we know that there exist natural numbers $q$ and $\bar{q}$ such that $2L + 1 = qL' + \bar{q}$, with $0 \le \bar{q} < q$. Note that this last equation can be written in the form

$$2L + 1 = \sum_{i=1}^{q} L'_i, \tag{49}$$

with $L'_i = L' + \theta_{\bar{q}}(i)$ and $\theta_{\bar{q}}(i) = 1$ for $i \le \bar{q}$ and $\theta_{\bar{q}}(i) = 0$ for $i > \bar{q}$. For our purposes, the relevant fact is that $L' \le L'_i \le L' + 1$. In the sequel, for any given pair of real numbers $a, b$ we will use the notation $[a, b]_l$ for $[a, b] \cap \mathbf{Z}$. We now will subdivide the box $[-L, L]_l$ in intervals according to (49). Thus, we define $I_1 := [-L, -L + L'_1 - 1]_l$ and for $1 < i \le q$ we let $I_i := [-L + \sum_{j=1}^{i-1} L'_i, -L + \sum_{j=1}^{i} L'_i - 1]_l$. Note that $I_q = [L - L'_q + 1, L]_l$ and $|I_i| = L'_i$. Next, we introduce a second parameter $r$ which is a natural number smaller than or equal to $L'$. We will call $r$ the *fine scale*. Then, for each $I_i$ we define an interval $J_i$ such that $J_i \subset I_i$, $|J_i| = L' - 2r$ and the endpoints of $J_i$ are at a distance larger than $r$ to the endpoints of $I_i$. To do so, first let $r_i := r + \theta_{\bar{q}}(i)$. Then define $J_1 := [-L + r, -L + L'_1 - 1 - r_1]_l$ and for $1 < i \le q$ we let $J_i := [-L + \sum_{j=1}^{i-1} L'_i + r, -L + \sum_{j=1}^{i} L'_i - 1 - r_i]_l$.

We now proceed to define the partition in $\Lambda_L$ in partition boxes and define the corresponding decomposition of the index set. First we define the set $\mathcal{I} := \{1, 2, \ldots, q\}^d$, which will correspond to the indexes parameterizing the sub-boxes. For a given element $\mathbf{i} \in \mathcal{I}$, of the form $\mathbf{i} = (i_1, \ldots, i_d)$ with $1 \le i_k \le q$, $1 \le k \le d$, we define

$$\Lambda'_{\mathbf{i}} := I_{i_1} \times I_{i_2} \times \cdots \times I_{i_d}.$$

We will call such a set a *partition box*. By definition the cardinality $|\Lambda'_{\mathbf{i}}|$ of a partition box satisfies

$$(L')^d \le |\Lambda'_{\mathbf{i}}| \le (L'+1)^d. \tag{50}$$

Note also that the partition boxes define a partition of $\Lambda_L$ so that $\Lambda_L = \bigcup_{\mathbf{i}\in\mathcal{I}} \Lambda'_{\mathbf{i}}$ where the union is disjoint.

Next we define a partition of the index set $\mathcal{I}$. Consider the collection $\mathcal{K}$ of subsets of $\{1, 2, \ldots, d\}$. Note that $|\mathcal{K}| = 2^d$. Now given $K \in \mathcal{K}$ we define $\mathcal{I}_K$ as the subset of $\mathcal{I}$ having coordinates which are even for $k \in K$ and odd for $k \notin K$. In other words, if we define $\mathbb{E}$ as the set of even natural numbers and $\mathbb{O}$ as the set of odd natural numbers, then

$$\mathcal{I}_K := \{\mathbf{i} = (i_1, \ldots, i_d) \in \mathcal{I} : i_k \in \mathbb{E} \text{ if } k \in K,\ i_k \in \mathbb{O} \text{ if } k \notin K,\ 1 \le k \le d\}.$$

Note that $\{\mathcal{I}_K : K \in \mathcal{K}\}$ defines a partition of $\mathcal{I}$ so that $\mathcal{I} = \bigcup_{K\in\mathcal{K}} \mathcal{I}_K$ is a disjoint union. Hence the union $\bigcup_{K\in\mathcal{K}} \bigcup_{\mathbf{i}\in\mathcal{I}_K} \Lambda'_{\mathbf{i}}$ is a disjoint union. We will refer to such a decomposition as the *parity partition at scale* $L'$ of $\Lambda_L$. Furthermore, given $K \in \mathcal{K}$ and any pair of boxes $\Lambda'_{\mathbf{i}}$ and $\Lambda'_{\mathbf{j}}$ with $\mathbf{i}, \mathbf{j} \in \mathcal{I}_K$ and $\mathbf{i} \ne \mathbf{j}$, we have that

$$\operatorname{dist}(\Lambda'_{\mathbf{i}}, \Lambda'_{\mathbf{j}}) \ge L'. \tag{51}$$

Here for any pair of subsets $A, B \subset \mathbb{Z}^d$ we define $\operatorname{dist}(A, B) := \inf_{x\in A, y\in B} |x - y|$. In other words (51) expresses the fact that the distance between any pair of partition boxes with different subindexes in $\mathcal{I}_K$ is larger than or equal to $L'$. This completes the description of the partition of $\Lambda_L$ in partition boxes.

Next, we describe the partition of $\Lambda_L$ in the strip set and main boxes. Given an $\mathbf{i} \in \mathcal{I}$ we let

$$\Lambda''_{\mathbf{i}} := J_{i_1} \times J_{i_2} \times \cdots \times J_{i_d}.$$

Such a box will be called a *main box*. Its cardinality is $|\Lambda''_{\mathbf{i}}| = (L' - 2r)^d$. Now let

$$S_L := \Lambda_L - \bigcup_{\mathbf{i}\in\mathcal{I}} \Lambda''_{\mathbf{i}}.$$

Such a set will be called the *strip set*. Note that $S_L$ and $\{\Lambda''_L : \mathbf{i} \in \mathcal{I}\}$ define a partition of $\Lambda_L$. We will refer to such a partition as the *strip-box partition at scale* $L'$ of $\Lambda_L$. We furthermore remark the following cardinality estimate for the strip set which will be useful later:

$$\frac{|S_L|}{(2L+1)^d} \le \frac{(L'+1)^d - (L'-2r)^d}{(L')^d}, \tag{52}$$

where we have used the fact that $|\mathcal{I}| = q^d$.

5.2. *Moment inequalities and decoupling.* Several inequalities and technical results obtained via the partition analysis will be used repeatedly throughout. Here we will derive some of these facts and recall some standard ones.

First we recall the following well-known inequality due to von Bahr and Esseen (see page 82, Exercise 2.6.20, of [8]).

LEMMA 5 (von Bahr–Esseen). *Let $X_1, \dots, X_n$ be mean-zero independent random variables and let $S_n := \sum_{i=1}^n X_i$. Then if $E$ denotes the expectation with respect to the joint law of $X_1, \dots, X_n$, and $1 \le r \le 2$, it is true that*

$$E|S_n|^r \le 2 \sum_{k=1}^n E|X_k|^r. \tag{53}$$

The following lemma is a corollary of (53) and (19).

LEMMA 6. *Consider the field of truncated quenched survival probabilities $\{\tilde{p}_1(x,t,w) : x \in \mathbb{Z}^d\}$ at scale 1. Let $L(t), L'(t) : [0,\infty) \to \mathbb{N}$ be functions such that $t \le L'(t) \le L(t)$. Then, if $1 \le r \le 2$, it is true that*

$$\begin{aligned} &\left\langle \left| \sum_{x \in \Lambda_L} (\tilde{p}_1 - \langle \tilde{p}_1 \rangle) \right|^r \right\rangle \\ &\qquad \le 2(2L'+2)^{d(r-1)}(2L+1)^d \langle |\tilde{p}_1 - \langle \tilde{p}_1 \rangle|^r \rangle. \end{aligned} \tag{54}$$

PROOF. We first perform the parity partition at scale $L'$ of the box $\Lambda_L$ so that $\Lambda_L = \bigcup_{K \in \mathcal{K}} \bigcup_{\mathbf{i} \in \mathcal{I}_{\mathcal{K}}} \Lambda_{\mathbf{i}}$. We will use the notation $c_x$ instead of $\tilde{p}_1(x,t,w) - \langle \tilde{p}_1 \rangle(x,t)$. Then,

$$\left\langle \left| \sum_{x \in \Lambda_L} c_x \right|^r \right\rangle = \left\langle \left| \sum_{K \in \mathcal{K}} \sum_{\mathbf{i} \in \mathcal{I}_K} \sum_{x \in \Lambda'_{\mathbf{i}}} c_x \right|^r \right\rangle. \tag{55}$$

Applying Jensen's inequality (19) to the right-hand side of (55) we see that $\langle | \sum_{K \in \mathcal{K}} \sum_{\mathbf{i} \in \mathcal{I}_K} \sum_{x \in \Lambda'_{\mathbf{i}}} c_x|^r \rangle \le 2^{d(r-1)} \sum_{K \in \mathcal{K}} \langle | \sum_{\mathbf{i} \in \mathcal{I}_K} \sum_{x \in \Lambda'_{\mathbf{i}}} c_x|^r \rangle$. Now, since the probabilities $\tilde{p}_1$ are truncated, by (51), and since $t \le L'(t) \le L(t)$, the random variables $\{\sum_{x \in \Lambda'_{\mathbf{i}}} c_x : \mathbf{i} \in \mathcal{I}_K\}$ are independent for every $K \in \mathcal{K}$. Therefore by the von Bahr–Esseen inequality (53), Jensen's inequality (19) and the volume estimate (50) we have that $2^{d(r-1)} \sum_{K \in \mathcal{K}} \langle | \sum_{\mathbf{i} \in \mathcal{I}_K} \sum_{x \in \Lambda'_{\mathbf{i}}} c_x|^r \rangle \le 2 2^{d(r-1)} \sum_{K \in \mathcal{K}} \sum_{\mathbf{i} \in \mathcal{I}_K} |\Lambda'_{\mathbf{i}}|^{r-1} \sum_{x \in \Lambda'_{\mathbf{i}}} \langle |c_x|^r \rangle \le 2(2L'+2)^{d(r-1)}(2L+1)^d \langle |c_0|^r \rangle$. Combining with (55), this completes the proof. □

We end this subsection with the following lemma which is an easy consequence of Lemma 1 and the asymptotic decay of the annealed survival probabilities.

LEMMA 7. *Consider the field of quenched survival probabilities* $\{p(x,t,w): x \in \mathbb{Z}^d\}$ *and that of truncated quenched survival probabilities* $\{\tilde p_1(x,t,w): x \in \mathbb{Z}^d\}$ *at scale* 1. *Let* $L(t): [0,\infty) \to \mathbb{N}$ *be such that* $L \gg 1$. *Then the following statements are true.*

(i) *For every* $t \geq 0$,

$$\left|\frac{p^L}{\langle p\rangle} - \frac{\sum_{x\in\Lambda_{L(t)}} \tilde p_1}{(2L+1)^d \langle \tilde p_1\rangle}\right| \leq \frac{k_1(1,d)e^{-k_2(1,d)t}}{\langle p\rangle}\left(1 + \frac{\sum_{x\in\Lambda_{L(t)}} \tilde p_1}{(2L+1)^d\langle\tilde p_1\rangle}\right). \tag{56}$$

(ii) *For every function* $f(t): [0,\infty) \to (0,\infty)$ *such that* $\log f(t) \gg -t^{d/(d+2)}$, *we have that*

$$\left|\frac{\sum_{x\in\Lambda_L}(p - \langle p\rangle)}{(2L+1)^d f(t)} - \frac{\sum_{x\in\Lambda_L}(\tilde p_1 - \langle\tilde p_1\rangle)}{(2L+1)^d f(t)}\right| \ll 1. \tag{57}$$

PROOF OF LEMMA 7(i). This is a direct consequence of the following inequality:

$$\begin{aligned}&\left|\frac{p^L}{\langle p\rangle} - \frac{\sum_{x\in\Lambda_{L(t)}}\tilde p_1}{(2L+1)^d\langle\tilde p_1\rangle}\right| \\ &\qquad \leq \frac{1}{(2L+1)^d\langle p\rangle\langle\tilde p_1\rangle}\left|\sum_{x\in\Lambda_L}(p - \tilde p_1)\right| + \frac{|\langle p_1 - p\rangle| \, |\sum_{x\in\Lambda_L}\tilde p_1|}{(2L+1)^d\langle p\rangle\langle\tilde p_1\rangle}\end{aligned}$$

and (13). □

PROOF OF LEMMA 7(ii). This follows again directly from (13). □

5.3. *The annealed asymptotics.* We proceed to prove the law of large numbers stated in (1). In other words we will prove that whenever $L(t) \geq \exp(\gamma\frac{c_2}{d}t^{d/(d+2)})$ for some $\gamma > \gamma_1$, where $\gamma_1 = \frac{2}{d+2}$, then in $\mu$-probability it is true that

$$\frac{p^{L(t)}(0,t,w)}{\langle p(0,t)\rangle} \sim 1. \tag{58}$$

To do so we first remark that by (56) and (21) it is enough to show that in $\mu$-probability,

$$\frac{\sum_{x\in\Lambda_L}\tilde p_1(x,t,w)}{(2L+1)^d\langle\tilde p_1(0,t)\rangle} \sim 1. \tag{59}$$

We will show that there is an $\varepsilon > 0$ such that

$$\left\langle\left|\frac{\sum_{x\in\Lambda_L}\tilde p_1}{(2L+1)^d\langle\tilde p_1\rangle} - 1\right|^{1+\varepsilon}\right\rangle \ll 1. \tag{60}$$

Remark that the right-hand side of (60) can be rewritten as

$$\left\langle\left|\frac{\sum_{x\in\Lambda_L}\tilde{p}_1}{(2L+1)^d\langle\tilde{p}_1\rangle}-1\right|^{1+\varepsilon}\right\rangle=\frac{\langle|\sum_{x\in\Lambda_L}(\tilde{p}_1-\langle\tilde{p}_1\rangle)|^{1+\varepsilon}\rangle}{(2L+1)^{d(1+\varepsilon)}\langle\tilde{p}_1\rangle^{1+\varepsilon}}. \tag{61}$$

At this point we make use of the parity partition decomposition for $\Lambda_L$ previously defined to deal with the numerator of the right-hand side of (61) via inequality (54) with $r=1+\varepsilon$. We will choose a time-dependent mesoscopic scale $L'(t):=\exp\{\gamma'\frac{c_2}{d}t^{d/(d+2)}\}$, where $0<\gamma'<\gamma-\gamma_1$. Therefore, the right-hand side of (61) is upper-bounded by

$$\frac{2(2L'+2)^{d\varepsilon}(2L+1)^d\langle|\tilde{p}_1-\langle\tilde{p}_1\rangle|^{1+\varepsilon}\rangle}{(2L+1)^{d(1+\varepsilon)}\langle\tilde{p}_1\rangle^{1+\varepsilon}}. \tag{62}$$

Now, by Lemma 4(i), we know that for $t\geq 0$ one has that

$$\langle\tilde{p}_1\rangle^{1+\varepsilon}=\exp\{-(1+\varepsilon)\times c_2t^{d/(d+2)}+o(t^{d/(d+2)})\}$$

(this particular case is a classical result of [5]). By part (ii) of the same corollary we also have that $\langle|\tilde{p}_1-\langle\tilde{p}_1\rangle|^{1+\varepsilon}\rangle=\exp\{-(1+\varepsilon)^{d/(d+2)}c_2t^{d/(d+2)}+o(t^{d/(d+2)})\}$. Hence, using these facts we see that (62) is upper-bounded by the expression $\frac{2^{d\varepsilon+1}(L'+1)^{d\varepsilon}}{(2L+1)^{d\varepsilon}}\exp\{[(1+\varepsilon)-(1+\varepsilon)^{d/(d+2)}]c_2t^{d/(d+2)}+o(t^{d/(d+2)})\}$. Then, from the inequality $(L'+1)^{d\varepsilon}/(2L+1)^{d\varepsilon}\leq\exp\{-\varepsilon(\gamma-\gamma')c_2t^{d/(d+2)}\}$ we see that (62) is upper-bounded by

$$\begin{aligned}&2^{d\varepsilon+1}\exp\{-[\varepsilon(\gamma-\gamma')\\&\qquad-((1+\varepsilon)-(1+\varepsilon)^{d/(d+2)})]c_2t^{d/(d+2)}+o(t^{d/(d+2)})\}.\end{aligned} \tag{63}$$

But since $\lim_{\varepsilon\to 0}\frac{(1+\varepsilon)-(1+\varepsilon)^{d/(d+2)}}{\varepsilon}=\frac{2}{d+2}=\gamma_1$, it follows that for $\varepsilon$ small enough, since $\gamma'<\gamma-\gamma_1$, the exponent $\varepsilon(\gamma-\gamma')-((1+\varepsilon)-(1+\varepsilon)^{d/(d+2)})$ is positive for $\gamma>\gamma_1$. Thus, we can choose $\varepsilon$ small enough to upper-bound (63) by $\exp\{-\mathrm{const}_1\,t^{d/(d+2)}+o(t^{d/(d+2)})\}$, where $\mathrm{const}_1>0$ is a constant depending only on $d$ and $\gamma$. This proves the validity of (60).

5.4. *The Gaussian asymptotics.* Here we will prove the central limit theorem stated in (4). This time we will need to perform a strip-box partition of the box $\Lambda_L$ into the strip set $S_L$ and the main boxes. We will choose the mesoscopic scale $L'(t)=\exp\{\gamma'\frac{c_2}{d}t^{d/(d+2)}\}$ with $0<\gamma'<\gamma$ and the fine scale $r=t^d$. Subsequently, $\gamma'$ will be chosen small enough. First note that by part (iv) of Lemma 3 and part (iv) of Lemma 4, both applied to the collection of sets $U_t=\Lambda_{L(t)}$, and Lemma 7 it is enough to prove that

$$\frac{\sum_{x\in\Lambda_L}(\tilde{p}_1(x,t)-\langle\tilde{p}_1(x,t)\rangle)}{\sqrt{\mathrm{Var}_\mu\sum_{x\in\Lambda_L}\tilde{p}_1(x,t)}} \tag{64}$$

converges in distribution to the normal law $\mathcal{N}(0,1)$. On the other hand,

$$\frac{\sum_{x\in\Lambda_L}(\tilde p_1-\langle\tilde p_1\rangle)}{\sqrt{\operatorname{Var}_\mu\sum_{x\in\Lambda_L}\tilde p_1}}=\frac{\sum_{x\in S_L}(\tilde p_1-\langle\tilde p_1\rangle)}{\sqrt{\operatorname{Var}_\mu\sum_{x\in\Lambda_L}\tilde p_1}}+\frac{\sum_{\mathbf{i}\in\mathcal{I}}\sum_{x\in\Lambda''_{\mathbf{i}}}(\tilde p_1-\langle\tilde p_1\rangle)}{\sqrt{\operatorname{Var}_\mu\sum_{x\in\Lambda_L}\tilde p_1}}. \tag{65}$$

We begin by showing that the strip component of the decomposition (65) converges to 0 in $\mu$-probability. To do so it is enough to prove that the variance of such a term converges to 0. But for $t$ large enough we have

$$\frac{\operatorname{Var}_\mu(\sum_{x\in S_L}\tilde p_1)}{\operatorname{Var}_\mu(\sum_{x\in\Lambda_L}\tilde p_1)}\sim\frac{|S_L|}{(2L+1)^d}\le 2d(2t)^d\exp\left(-\gamma'\frac{c_2}{d}t^{d/(d+2)}\right) \tag{66}$$

where we have used first, part (iv) of Lemma 4 with $U_t=S_{L(t)}$ and $U_t=\Lambda_{L(t)}$, and in the inequality we have used estimate (52).

We have thus reduced the proof to showing that the second term of the decomposition of the right-hand side of inequality (65) converges in distribution to $\mathcal{N}(0,1)$. By the choice of the fine scale $r=t^d$, the random variables $\{\sum_{x\in\Lambda''_{\mathbf{i}}}(\tilde p_1-\langle\tilde p_1\rangle):\mathbf{i}\in\mathcal{I}\}$ are independent, so that considering again the estimate (66), it is enough to verify the following version of the Lyapunov condition. There is an $\varepsilon>0$ such that

$$\frac{\sum_{\mathbf{i}\in\mathcal{I}}\langle|\sum_{x\in\Lambda''_{\mathbf{i}}}(\tilde p_1-\langle\tilde p_1\rangle)|^{2+\varepsilon}\rangle}{(\sum_{\mathbf{i}\in\mathcal{I}}\operatorname{Var}_\mu\sum_{x\in\Lambda''_{\mathbf{i}}}\tilde p_1)^{1+\varepsilon/2}}\ll 1. \tag{67}$$

Now, by Jensen's inequality (19) with $r=2+\varepsilon$, the fact that $|\Lambda''_{\mathbf{i}}|\le|\Lambda'_{\mathbf{i}}|$, inequality (50), and the fact that $\langle|\tilde p_1-\langle\tilde p_1\rangle|^{2+\varepsilon}\rangle=\exp\{-(2+\varepsilon)^{d/(d+2)}c_2t^{d/(d+2)}+o(t^{d/(d+2)})\}$ [which follows from part (ii) of Lemma 4], we see that the numerator of the left-hand side of (67) is upper-bounded by $\exp\{-(2+\varepsilon)^{d/(d+2)}c_2t^{d/(d+2)}+o(t^{d/(d+2)})\}(L'+1)^{d(1+\varepsilon)}(2L+1)^d$. Furthermore, by part (iv) of Lemma 4 applied for each $\mathbf{i}\in\mathcal{I}$ to the collection of sets $U_t=\Lambda''_{\mathbf{i}}$, the quantity $\exp\{-2^{d/(d+2)}(1+\frac{\varepsilon}{2})c_2t^{d/(d+2)}+o(t^{d/(d+2)})\}(2L+1)^{d(1+\varepsilon/2)}$ divided by the denominator of the left-hand side of (67) converges to 1. Hence the left-hand side of (67) is upper bounded by

$$\begin{aligned}\exp\bigg\{\bigg((1+\varepsilon)\gamma'-\varepsilon\frac{\gamma}{2}+\bigg[2^{d/(d+2)}\Big(1+\frac{\varepsilon}{2}\Big)\\ -(2+\varepsilon)^{d/(d+2)}\bigg]\bigg)c_2t^{d/(d+2)}+o(t^{d/(d+2)})\bigg\}.\end{aligned} \tag{68}$$

Now, considering that $\lim_{\varepsilon\to0}\frac{2^{d/(d+2)}(1+\varepsilon/2)-(2+\varepsilon)^{d/(d+2)}}{\varepsilon}=2^{d/(d+2)}\frac{1}{2}\frac{2}{d+2}=\frac{1}{2}\gamma_2$, we see that for $\varepsilon$ and $\gamma'$ small enough and $\gamma>\gamma_2$, the exponent $(1+\varepsilon)\gamma'-\varepsilon\frac{\gamma}{2}+[2^{d/(d+2)}(1+\frac{\varepsilon}{2})-(2+\varepsilon)^{d/(d+2)}]$ of the bound (68) is negative. This proves the Lyapunov condition (67).

5.5. *The critical asymptotics.* In this subsection we will prove part (i) of Theorem 2. Let us note that as a corollary of (8) we obtain (3), showing the absence of an annealed behavior, for $L(t) \le \exp\{\gamma \frac{c_2}{d} t^{d/(d+2)}\}$ with $\gamma \le \gamma_1$. In fact, statement (8) implies (3), since $\langle p\rangle \sim \exp\{-c_2 t^{d/(d+2)} + o(t^{d/(d+2)})\}$ while $\exp\{-c_2 t^{d/(d+2)} + o(t^{d/(d+2)})\} \ge \exp\{-a(\gamma) c_2 t^{d/(d+2)}\}$.

Let us proceed with the proof of (8). Here we will not make use of the partition analysis. Note that it is enough to show that for every $0 < \delta < 1$ and $0 < \gamma \le \gamma_1$ one has that $\mu$-a.s.,

$$\lim_{n\to\infty} \sup_{n\le s\le n+1} \frac{\sum_{x\in\Lambda_{L(s)}} p(x,s)}{(2L(s)+1)^d \exp\{-(a(\gamma)-\delta)c_2 s^{d/(d+2)}\}} = 0. \tag{69}$$

Let $b = (\frac{1}{\alpha'}(a-\gamma))^{-\alpha/\alpha'} > 0$. Note also that $b = (\frac{d+2}{2}\gamma)^{(d+2)/d} \le 1$ whenever $0 \le \gamma \le \gamma_1$. By the Borel–Cantelli lemma, to prove (69) it suffices to prove that

$$\begin{aligned} &\left\langle \left( \sup_{n\le s\le n+1} \frac{\sum_{x\in\Lambda_{L(s)}} p(x,s)}{(2L(s)+1)^d \exp\{-(a(\gamma)-\delta)c_2 s^{d/(d+2)}\}} \right)^b \right\rangle \\ &\qquad \ll \exp\left\{ -\frac{b\delta c_2 n^{d/(d+2)}}{2} \right\}. \end{aligned} \tag{70}$$

Now, by the monotonicity of $p(x,s)$ and by (20) with $r = b$ applied in the summation of (70) we conclude that

$$\begin{aligned} &\left\langle \left( \sup_{n\le s\le n+1} \frac{\sum_{x\in\Lambda_L} p(x,s)}{(2L(s)+1)^d e^{-(a(\gamma)-\delta)c_2 s^{d/(d+2)}}} \right)^b \right\rangle \\ &\qquad \le \frac{(2L(n+1)+1)^d \langle p(0,n)^b\rangle}{(2L(n)+1)^{db} e^{-b(a(\gamma)-\delta)c_2 (n+1)^{d/(d+2)}}}. \end{aligned} \tag{71}$$

By Lemma 3(i) with $r = b$, we know that $\langle p(0,n)^b\rangle \le \exp\{-c_2 b^{d/(d+2)} n^{d/(d+2)} + o(n^{d/(d+2)})\}$. Thus, the right-hand side of (71) is upper-bounded by

$$\exp\{-[b^{d/(d+2)} - b(a(\gamma)-\delta) - \gamma(1-b)]c_2 n^{d/(d+2)} + o(n^{d/(d+2)})\}. \tag{72}$$

Here, we will analyze briefly the function $f(x) := x^{d/(d+2)} - xc - \gamma(1-x)$, appearing in the exponent of the exponential of (72) with $x = b$ and $c = a(\gamma) - \delta$. Note that $f(x)$ has no roots and is strictly negative for $c > a(\gamma)$, while it has one root for $c = a(\gamma)$ and two roots for $c < a(\gamma)$. Therefore, if $c = a(\gamma) - \delta$ we can find a value of $x$ for which $f(x)$ is positive, for example, $x = b$. Note also that $f(b) = b\delta > 0$ and the expression of display (72) is equal to $\exp\{-b\delta c_2 n^{d/(d+2)} + o(n^{d/(d+2)})\}$. This proves (70).

5.6. *The annealed non-Gaussian asymptotics.* Here we will prove the statement of (9). As a corollary we will obtain (5). To see this, let us observe that $\log\sqrt{\mathrm{Var}_\mu(p^{L(t)}(0,t,w))} \sim -(2^{-2/(d+2)} + \frac{\gamma}{2})c_2 t^{d/(d+2)}$. On the other hand, note that $\inf_{\gamma>0}(a(\gamma) - \frac{\gamma}{2}) = 2^{-2/(d+2)}$ and the infimum is attained at $\gamma = \gamma_2 = 2^{d/(d+2)}\frac{2}{d+2}$. Hence, $a(\gamma) \geq 2^{-2/(d+2)} + \frac{\gamma}{2}$. Thus, $\sqrt{\mathrm{Var}_\mu(p^{L(t)}(0,t,w))} \geq \exp\{-a(\gamma)c_2 t^{d/(d+2)}\}$. First note that by Lemma 7, we have

$$\frac{p^{L(t)} - \langle p\rangle}{e^{-(a(\gamma)-\delta)c_2 t^{d/(d+2)}}} \sim \frac{\sum_{x\in\Lambda_L}(\tilde{p}_1 - \langle\tilde{p}_1\rangle)}{(2L+1)^d e^{-(a(\gamma)-\delta)c_2 t^{d/(d+2)}}}. \tag{73}$$

So it is enough to show that the right-hand side of (73) converges to 0 in probability. Let $b = (\frac{1}{\alpha'}(a-\gamma))^{-\alpha/\alpha'}$. Note that $1 < b = (\frac{d+2}{2}\gamma)^{(d+2)/d} \leq 2$ for $\gamma_1 < \gamma < \gamma_2$. We will show that

$$\left\langle\left|\frac{\sum_{x\in\Lambda_L}(\tilde{p}_1 - \langle\tilde{p}_1\rangle)}{(2L+1)^d \exp\{-(a(\gamma)-\delta)c_2 t^{d/(d+2)}\}}\right|^b\right\rangle \ll 1. \tag{74}$$

To estimate the left-hand side of (74) we will make use of the parity partition of $\Lambda_L$. We will choose the mesoscopic scale $L' := \exp\{-\frac{\gamma'}{d}c_2 t^{d/(d+2)}\}$. Now, by (54) with $r = b$ we see that $(\sum_{x\in\Lambda_L}(\tilde{p}_1 - \langle\tilde{p}_1\rangle))^b \leq 2(2L'+2)^{d(b-1)}(2L+1)^d\langle|\tilde{p}_1 - \langle\tilde{p}_1\rangle|^b\rangle$. On the other hand, part (iii) of Lemma 4 tells us that $\langle|\tilde{p}_1 - \langle\tilde{p}_1\rangle|^b\rangle = \exp\{-b^{d/(d+2)}c_2 t^{d/(d+2)}\} + o(t^{d/(d+2)})$. We can thus bound the left-hand side of (74) by

$$\begin{aligned}&\exp\{-[b^{d/(d+2)} - b(a(\gamma)-\delta)\\ &\qquad - \gamma(1-b) - \gamma'(b-1)]c_2 t^{d/(d+2)} + o(t^{d/(d+2)})\}.\end{aligned} \tag{75}$$

As in the previous subsection we obtain the function $f(x) := x^{d/(d+2)} - xc - \gamma(1-x)$ in the exponent of the exponential of (75) with $x = b$ and $c = a(\gamma) - \delta$. The same analysis shows us that $f(b) = b\delta > 0$ and that the expression of (75) equals

$$\exp\{-(b\delta - \gamma'(b-1))c_2 t^{d/(d+2)} + o(t^{d/(d+2)})\}.$$

Choosing $\gamma' < \delta$, it is clear that the above expression converges to 0.

**6. The one-dimensional convergence to infinitely divisible laws.** In this section we will establish Theorem 3, describing the convergence in distribution to infinitely divisible laws of the quantities $\frac{p^L(0,t,w)}{\langle p^L(0,t)\rangle}$ and $\frac{p^L(0,t,w)-\langle p^L(0,t)\rangle}{\mathrm{Var}_\mu(p^L(0,t))}$ in the critical nonannealed and non-Gaussian asymptotics. The proof of this theorem will be given in two steps: a reduction to a problem of sums of independent random variables; and a straightforward application of classical results for the convergence of sums of independent random variables to infinitely divisible laws. The verification of the conditions giving convergence will be similar to the results of [3], where positive exponentials of random variables with Weibull-type and Fréchet-type tails were studied.

6.1. *Reduction to sums of independent random variables.* We first need to introduce some minimal notation. Given an environment $w \in X$, recall the definition of the obstacle set $\mathcal{G}(w) = \{y \in \mathbb{Z} : w(y) = 1\}$. We want to make use of the natural order of $\mathbb{Z}$ to enumerate the elements of $\mathcal{G}(w)$. Let us define $y_0 := \inf\{y \in \mathcal{G}(w) : y \geq 0\}$, the site having an obstacle with nonnegative but smallest coordinate. Then recursively in $m \in \mathbb{N}$ define $y_m := \inf\{y \in \mathcal{G}(w) : y > y_{m-1}\}$ and $y_{-m} := \sup\{y \in \mathcal{G}(w) : y < y_{-(m-1)}\}$. Note that for every integer $m$ one has that $y_m < y_{m+1}$ and furthermore $\mathcal{G}(w) = \bigcup_{i \in \mathbb{Z}}\{y_m\}$. Let us also define for each $m \in \mathbb{Z}$ the intervals $I_m(w) := [y_{m-1} + 0.5, y_m - 0.5]_l$. We denote their lengths $l_m := |I_m| = y_m - y_{m-1} - 1$. Note that $\{l_m(w) : m \in \mathbb{Z}/\{0\}\}$ is a set of independent identically distributed random variables with a geometric law of parameter $1 - p$ so that $\mu(l_m = k) = p(1-p)^k$ for every $m \in \mathbb{Z}$ and $k \geq 0$. On the other hand, $l_0$ is independent from the set $\{l_m(w) : m \in \mathbb{Z}/\{0\}\}$ but has a law given by $l_0 = l_0^+ + l_0^-$, where $l_0^+ := |I_0 \cap [0, \infty)|$ and $l_0^- := |I_0 \cap (-\infty, -1]|$ are independent identically distributed random variables with a geometric law of parameter $1 - p$. Also, at this point we recall the notation of Section 3, so that for each $m \in \mathbb{Z}$, $\{\lambda_n(I_m) : 0 \leq n \leq l_m\}$ and $\{\psi_n^{I_m} : 0 \leq n \leq l_m\}$ represent the Dirichlet eigenvalues and eigenfunctions, respectively, of the discrete Laplacian operator (15) on $I_m$. Here we have used the fact that there are $l_m$ different eigenvalues each one of geometric multiplicity 1. Our objective here is to prove the following proposition.

PROPOSITION 1. *Let $L(t) : [0, \infty) \to \mathbb{N}$ be such that $L(t) \geq 1$ and assume that $0 \leq p < 1$. Consider the field of quenched survival probabilities $\{p(x, t, w) : x \in \mathbb{Z}\}$ in dimension $d = 1$. Then for every $\varepsilon > 0$ $\mu$-a.s. eventually in $t$ it is true that*

$$\sum_{k=-m}^{m} X_k \leq \sum_{x \in \Lambda_{L(t)}} p(x, t, w) \leq \sum_{k=-M}^{M} X_k, \tag{76}$$

*where*

$$X_k := \frac{2(l_k + 1)}{\ell_1} \exp\left\{-4t \frac{\ell_1}{(l_k + 1)^2}(1 + o_1(l_k))\right\}(1 + o_4(l_k, t)), \tag{77}$$

*with $m := [(1-\varepsilon)L\frac{p}{1-p}]$, $M := [(1+\varepsilon)L\frac{p}{1-p}]$, $\ell_1 := \frac{\pi^2}{8}$ is the principal Dirichlet eigenvalue of the operator $-\frac{1}{2}\Delta$ on the ball of radius 1, $o_4(x, t) : [0, \infty)^2 \to \mathbb{R}$* *is a function such that* $|o_4(x,t)| \le \frac{200}{1+x} + 500e^{-t(2\pi^2/(x+1)^2)(1+o_2(x))}$ *and* $o_1(x), o_2(x) : [0,\infty) \to \mathbb{R}$ *are functions such that* $|o_1(x)| \le \frac{\pi^2}{12(1+x)^2}$ *and* $|o_2(x)| \le \frac{10}{(1+x)^2}$.

To prove Proposition 1 our first step will be the following lemma which tells us that the averaged survival probability at scale $L$ is essentially equal to the averaged survival probability at the random scale $y_{(1-p)L}$.

LEMMA 8. *Let* $L(t) : [0,\infty) \to \mathbb{N}$ *be such that* $L(t) \ge 1$ *and assume that* $0 \le p < 1$. *Consider the field of quenched survival probabilities* $\{p(x,t,w) : x \in \mathbb{Z}\}$ *in dimension* $d = 1$. *Then for every* $\varepsilon > 0$ $\mu$*-a.s. eventually in* $t$ *it is true that*

$$\sum_{x=y_{-m}}^{y_m} p(x,t,w) \le \sum_{x\in\Lambda_{L(t)}} p(x,t,w) \le \sum_{x=y_{-M}}^{y_M} p(x,t,w), \tag{78}$$

*where* $m := [(1-\varepsilon)L\frac{p}{1-p}]$ *and* $M := [(1+\varepsilon)L\frac{p}{1-p}]$.

PROOF. Note that for every nonnegative natural $n$, $y_n = l_0^+ + \sum_{i=1}^n (l_n + 1)$. Since $l_0^+, l_1, \dots$ is a sequence of independent geometric random variables of parameter $1-p$, by the strong law of large numbers since $L(t) \gg 1$ we know that for every $\varepsilon > 0$, $\mu$-a.s. eventually in $t$ we have that $y_m < L(t) < y_M$. Similarly, $\mu$-a.s. eventually in $t$, $y_{-M} < -L(t) < y_{-m}$. This proves (78). □

Next, we recall some elementary estimates for the principal Dirichlet eigenvalue and the $L^1$ norm of the principal Dirichlet eigenfunction on an interval in terms of the length of the interval.

LEMMA 9. *Let* $I \subset \mathbb{Z}$ *be a bounded nonempty interval. Consider the Dirichlet eigenvalues* $\{\lambda_n(I) : n \in \mathcal{U}\}$ *and eigenfunctions* $\{\psi_n^I : n \in \mathcal{U}\}$ *of the discrete Laplacian* $\Delta_d$ *on* $I$, *where* $\mathcal{U}$ *is a finite index set. Then:*

(i)

$$\lambda_0(I) = \frac{4\ell_1}{(l+1)^2}(1 + o_1(l)).$$

(ii)

$$\lambda_1(I) - \lambda_0(I) = \frac{12\ell_1}{(l+1)^2}(1 + o_2(l)).$$

(iii)

$$(\psi_0^I, \mathbf{1}_I) = \sqrt{\frac{2(l+1)}{\ell_1}}(1 + o_3(l)),$$

*where* $l=|I|$, $\ell_1=\frac{\pi^2}{8}$ *is the principal Dirichlet eigenvalue of the differential operator* $-\frac{1}{2}\Delta$ *on the ball of radius* 1 *and* $o_1(x),o_2(x),o_3(x):[0,\infty)\to[0,\infty)$ *are functions such that* $|o_1(x)|\le\frac{\pi^2}{12(1+x)^2}$, $|o_2(x)|\le\frac{10}{(1+x)^2}$ *and* $|o_3(x)|\le\frac{10}{1+x}$.

PROOF. By translation invariance, without loss of generality we can assume that $I=\{0,1,\dots,l-1\}$. Recall that $\lambda_n(I)=1-\cos(\frac{(n+1)\pi}{l+1})$ for $n\in\mathcal{U}$ with $\mathcal{U}=\{0,1,\dots,l-1\}$, and that $\psi_0^I(x)=\sqrt{\frac{2}{l+1}}\sin(\frac{\pi}{l+1}(x+1))$ for $x\in I$. Parts (i) and (ii) follow directly from the fact that $\frac{x^2}{2}(1-\frac{x^2}{12})\le 1-\cos x\le\frac{x^2}{2}$. To deduce part (iii) note that since $\Delta_1\psi_0^I(x)=\lambda_0\psi_0^I(x)$ when $x\in I$ [recall the definition (15) of $\Delta_1$], then $(\psi_0^I,\mathbf{1}_I)=\sum_{x=0}^{l-1}\psi_0^I(x)=\frac{\psi_0^I(0)+\psi_0^I(l-1)}{2\lambda_0(I)}$. Thus,

$$(\psi_0^I,\mathbf{1}_I)=\sqrt{\frac{2}{l+1}}\left(\frac{\sin\pi/(l+1)}{1-\cos\pi/(l+1)}\right).$$

We can now deduce (iii) from the mentioned bounds for $\cos x$ and the inequalities, $x-\frac{1}{6}x^3\le\sin x\le x$. $\square$

As a consequence of the estimates of Lemma 9 we deduce the following result, which together with Lemma 8 finishes the proof of Proposition 1.

LEMMA 10. *Consider the field of quenched survival probabilities* $\{p(x,t,w):x\in\mathbb{Z}\}$ *in dimension* $d=1$. *Then for every* $r\in\mathbb{N}$, *we have that*

$$\sum_{x=y_{-r}}^{y_r}p(x,t,w)=\sum_{k=-r}^{r}\frac{2(l_k+1)\mathbf{1}_{\mathbb{N}^+}(l_k)}{\ell_1}\times\exp\left\{-t\frac{4\ell_1}{(l_k+1)^2}(1+o_1(l_k))\right\}(1+o_4(l_k,t)),\tag{79}$$

*where* $\mathbf{1}_{\mathbb{N}^+}$ *is the indicator function of the set of natural numbers larger than* 0, $\ell_1:=\frac{\pi^2}{8}$, $o_4(x,t):[0,\infty)^2\to\mathbb{R}$, $|o_4(x,t)|\le\frac{200}{1+x}+500e^{-t(2\pi^2/(x+1)^2)(1+o_2(x))}$ *and* $o_1(x),o_2(x):[0,\infty)\to\mathbb{R}$ *are functions such that* $|o_1(x)|\le\frac{\pi^2}{12(1+x)^2}$ *and* $|o_2(x)|\le\frac{10}{(1+x)^2}$.

PROOF. Let $k$ be such that $-m\le k\le m$. Note that if there is an $x\in I_k$ (which implies that $l_k=|I_k|\ge 1$), we have the expansion

$$p(x,t,w)=\sum_{n=0}^{l_k-1}e^{-t\lambda_n(I_k)}(\psi_n^{I_k},\mathbf{1}_{I_k})\psi_n^{I_k}(x).$$

Hence, $\sum_{x\in I_k} p(x,t,w) = \sum_{n=0}^{l_k-1} e^{-t\lambda_n(I_k)}(\psi_n^{I_k}, \mathbf{1}_{I_k})^2$ and

$$\sum_{x\in I_k} p(x,t,w) = e^{-t\lambda_0(I_k)}\left((\psi_0^{I_k}, \mathbf{1}_{I_k})^2 + \sum_{n=1}^{l_k-1} e^{-t(\lambda_n(I_k)-\lambda_0(I_k))}(\psi_n^{I_k}, \mathbf{1}_{I_k})^2\right). \tag{80}$$

Now, the summation in the right-hand side of (80) is bounded above by $e^{-t(\lambda_1(I_k)-\lambda_0(I_k))}(l_k - (\psi_0^{I_k}, \mathbf{1}_{I_k})^2)$. Hence we conclude that

$$e^{-t\lambda_0(I_k)}A_k \le \sum_{x\in I_k} p(x,t,w) \le e^{-t\lambda_0(I_k)}(A_k + e^{-t(\lambda_1(I_k)-\lambda_0(I_k))}(l_k - A_k)),$$

where $A_k := (\psi_0^{I_k}, \mathbf{1}_{I_k})^2$. Summing up over $k$ and applying the estimates of Lemma 9, we obtain (79). □

6.2. *Convergence to infinitely divisible laws.* We now wish to describe how to finish the proof of Theorem 3 via Proposition 1. By Proposition 1 it is clear that it is enough to prove the following.

PROPOSITION 2. *Let* $L(t) := e^{\nu[(\gamma c_2/\nu)t^{1/3}]_-}$, $c > 0$ *and* $r := [cL\frac{p}{1-p}]$. *Consider a sequence* $\{l_k : k \in \mathbb{N}\}$ *of independent identically distributed geometric random variables of parameter* $1-p$ *and joint law* $\upsilon$. *Let* $X_k(t) := \frac{2(l_k+1)}{\ell_1}\exp\{-t\frac{4\ell_1}{(l_k+1)^2}(1+o_1(l_k))\}(1+o_4(l_k,t))$, *where* $o_1, o_2, o_4 : [0,\infty) \to \mathbb{R}$ *are functions such that* $|o_4(x,t)| \le \frac{200}{1+x} + 500e^{-t(2\pi^2/(x+1)^2)(1+o_2(x))}$, $|o_1(x)| \le \frac{\pi^2}{12(1+x)^2}$ *and* $|o_2(x)| \le \frac{10}{(1+x)^2}$. *Then the following statements are true.*

(i) *If* $\gamma < \gamma_1 := 1/3$, *then*

$$\lim_{t\to\infty} \frac{1}{s_1(\gamma)t^{1/3}} \exp\left\{\frac{4t\ell_1}{[\gamma c_2 t^{1/3}/\nu]_-^2}\right\} \sum_{k=-r}^{r} X_k(t) = X_{\beta_1,0,\mathcal{L}}.$$

(ii) *If* $\gamma_1 < \gamma < \gamma_2 := 2^{1/3}/3$, *then*

$$\lim_{t\to\infty} \frac{1}{s_1(\gamma)t^{1/3}} \exp\left\{\frac{4t\ell_1}{[\gamma c_2 t^{1/3}/\nu]_-^2}\right\} \sum_{k=-r}^{r} (X_k(t) - \upsilon(X_1(t))) = X_{\beta_2,0,\mathcal{L}}.$$

*In both cases the convergence is in distribution,* $s_1(\gamma) := \frac{\gamma c_2}{\ell_1 \nu}$, $\beta_1 := \frac{2cp^2}{(1-p)^2} \times \sum_{k=-\infty}^{\infty} \frac{(1-p)^k}{(1-p)^{k/a_1} + (1-p)^{-k/a_1}}$, $\beta_2 := \frac{2cp^2}{(1-p)^2}\sum_{k=-\infty}^{\infty} \frac{(1-p)^{k(1+2/a_1)}}{(1-p)^{k/a_1}+(1-p)^{-k/a_1}}$, *the Lévy spectral function* $\mathcal{L}(x) := -\frac{2cp}{1-p}(1-p)^{-[(1/\nu)\log x^{a_1(\gamma)}]}$ *for* $x > 0$ *and* $\mathcal{L}(x) = 0$ *for* $x < 0$ *and* $a_1 := (\frac{\alpha}{\alpha'}\gamma)^{1/\alpha'}$.

The proof of Proposition 2 will be the content of the following subsections. We will verify the conditions according to the classical results (see, e.g., Theorem 3.3 of [8]). Let us recall them.

THEOREM 4. *Let* $n(t):[0,\infty)\to\mathbb{N}$ *and for each* $t$ *let* $\{Y_k(t):1\le k\le n(t)\}$ *be a sequence of independent identically distributed random variables. Call* $P_t$ *the law of* $Y_1(t)$. *Assume that for every* $\varepsilon>0$ *it is true that*

$$\lim_{t\to\infty} P_t(Y_1(t)>\varepsilon)=0.$$

*Now let* $\mathcal{L}(x):\mathbb{R}/\{0\}\to\mathbb{R}$ *be a Lévy spectral function,* $\beta\in\mathbb{R}$ *and* $\sigma>0$. *Then the following statements are equivalent:*

(i)

$$\lim_{t\to\infty}\sum_{k=1}^{n(t)} Y_k(t)=X_{\beta,\sigma,\mathcal{L}},$$

*where the convergence is in distribution.*

(ii) *Define for* $\tau>0$ *the truncated random variable at level* $\tau$ *as* $Z_\tau(t):=Y_1(t)\mathbf{1}_{|Y_1(t)|\le\tau}$. *Also, let* $E_t(\cdot)$ *and* $\operatorname{Var}_t(\cdot)$ *denote the expectation and variance corresponding to the law* $P_t$. *Then if* $x$ *is a continuity point of* $\mathcal{L}$,

$$\mathcal{L}(x)=\begin{cases}\lim_{t\to\infty} n(t)P_t(Y_1(t)\le x), & \text{for } x<0,\\ -\lim_{t\to\infty} n(t)P_t(Y_1(t)>x), & \text{for } x>0,\end{cases}\tag{81}$$

$$\sigma^2=\lim_{\tau\to 0}\lim_{t\to\infty} n(t)\operatorname{Var}_t(Z_\tau(t)),\tag{82}$$

*and for any* $\tau>0$ *which is a continuity point of* $\mathcal{L}(x)$,

$$\begin{aligned}\beta=&\lim_{n\to\infty} n(t)E_t(Z_\tau(t))\\ &+\int_{|x|>\tau}\frac{x}{1+x^2}\,d\mathcal{L}(x)-\int_{\tau\ge|x|>0}\frac{x^3}{1+x^2}\,d\mathcal{L}(x).\end{aligned}\tag{83}$$

In the next subsections, we proceed to verify conditions (81), (82) and (83), for the triangular array given by $X_k(t)$ in part (i) of Proposition 2 and by $X_k(t)-v(X_k(t))$ in part (ii) of the same proposition.

6.3. *The Lévy spectral function in the critical case.* Here we begin the proof of part (i) of Proposition 2 by identifying the Lévy spectral function $\mathcal{L}(x)$ verifying condition (81) with

$$Y_1(t):=\frac{1}{s_1t^{1/3}}\exp\biggl\{\frac{4t\ell_1}{[\gamma c_2t^{1/3}/\nu]_-^2}\biggr\}X_1(t),$$

and $n(t) := 2r(t)+1$ with $r(t) := [cL(t)\frac{p}{1-p}]$, $L(t) := e^{\nu[(\gamma c_2/\nu)t^{1/3}]_-}$ and $c>0$.

Namely, we wish to prove that

$$(84)\qquad \begin{aligned} &-\lim_{t\to\infty} n(t)\upsilon(Y_1(t)>x) \\ &= \mathcal{L}(x) = \begin{cases} 0, & \text{for } x<0, \\ -\dfrac{2cp}{1-p}(1-p)^{-[(1/\nu)\log x^{a_1(\gamma)}]}, & \text{for } x>0. \end{cases} \end{aligned}$$

Since $X_1(t)\geq 0$ it is obvious that $\mathcal{L}(x)=0$ for $x<0$. Thus, we concentrate on the case $x>0$. Let us first note that if $l_1\leq t^{2/9}$, then $Y_1(t)\leq \exp\{-4\ell_1 t^{5/9}+o(t^{-5/9})\}$. It follows that $n(t)\upsilon(Y_1(t)>x, l_1\leq t^{2/9})=0$ for $t$ large enough. Hence,

$$(85)\qquad \lim_{t\to\infty} n(t)\upsilon(Y_1(t)>x) = \lim_{t\to\infty} n(t)\upsilon(Y_1(t)>x, l_1>t^{2/9}).$$

Let us also remark that $\upsilon(l_1\geq t^{2/5}) = \exp\{-\nu t^{2/5}+o(t^{-2/5})\}$. Thus, since $n(t)=2r(t)+1=\exp\{\gamma c_2 t^{1/3}+o(t^{-1/3})\}$ we have that $\lim_{t\to\infty} n(t)\upsilon(l_1\geq t^{2/5})=0$. Thus, from (85) we see that

$$(86)\qquad \lim_{t\to\infty} n(t)\upsilon(Y_1(t)>x) = \lim_{t\to\infty} n(t)\upsilon(Y_1(t)>x, t^{2/9}<l_1<t^{2/5}).$$

Now note that if $b_1$ and $b_2$ are positive real numbers, the function $f(y) = b_1 y e^{-b_2/y^2}$ is increasing. Therefore, if $y_0$ is the solution of $f(y)=x$ and $n$ is a natural number, we have that $f(n)>x$ if and only if $n>[y_0]$. Note that $y_0 = \sqrt{\frac{b_2}{\log b_1 y_0 - \log x}}$. This implies that if $b_1>x$, then either $y_0\leq 1$ or if $y_0>1$ then $y_0\leq \sqrt{\frac{b_2}{\log b_1-\log x}}$. Hence if $b_1>x$, then $y_0\leq \max\{\sqrt{\frac{b_2}{\log b_1-\log x}},1\}$. Thus,

$$(87)\qquad \begin{aligned} &\sqrt{\frac{b_2}{\log b_1 + (1/2)\log(b_2/(\log b_1-\log x)) - \log x}} \leq y_0 \\ &= \sqrt{\frac{b_2}{\log b_1 + 1/2\log(b_2(\log b_1 y_0 - \log x)) - \log x}}, \end{aligned}$$

whenever $b_1>x$ and $\sqrt{\frac{b_2}{\log b_1-\log x}}\geq 1$. Now, the inequality $Y_1(t)>x$ appearing in (86) can be expressed as $f(l_1+1)>x$ with $f(y)=b_1 y e^{-b_2/y^2}$, $b_1 := \frac{2}{\ell_1 s_1 t^{1/3}}\exp\{\frac{4t\ell_1}{[\gamma c_2 t^{1/3}/\nu]_-^2}\}(1+o_4(l_1,t))$ and $b_2 := 4t\ell_1(1+o_1(l_1))$. Hence,

$$(88)\qquad \upsilon(Y_1(t)>x, t^{2/9}<l_1<t^{2/5}) = \upsilon(l_1+1>[y_0], t^{2/9}<l_1<t^{2/5}).$$

But when $t^{2/9}<l_1<t^{2/5}$ we have that $|o_1(l_1)|\leq \frac{\pi^2}{12t^{4/9}}$ and $|o_4(l_1,t)|\leq \frac{200}{t^{2/9}} + 500e^{-\pi^2 t^{1/5}}$. Hence, $\log b_1 = (\frac{4t\ell_1}{[\gamma c_2 t^{1/3}/\nu]_-^2} - \log(s_1\ell_1 t^{1/3}/2)) + O(t^{-2/9})$, $b_2 =$

$4t\ell_1(1+O(t^{-4/9}))$ and the lower bound of (87) is of the order of $O(t^{1/3})$. Furthermore, for $t$ large enough we have $b_1 > x$ and $\sqrt{\frac{b_2}{\log b_1 - \log x}} \geq 1$ so that we can apply such an inequality. It follows from (87) and (88) that for $t$ large enough,

$$\upsilon(Y_1(t) > x, t^{2/9} < l_1 < t^{2/5}) = \upsilon(l_1 + 1 > [y_0]) + o(1/n(t)), \tag{89}$$

where we have used the fact that $\upsilon(l_1 \geq t^{2/5}) = o(1/n(t))$. Now, $= \upsilon(l_1 + 1 > [y_0]) = e^{-\nu[y_0]}$. We therefore need to get a good estimate on $y_0$. Furthermore, a short computation enables us to conclude that $\log \frac{b_2}{\log b_1 - \log x} = 2\log(\gamma c_2 t^{1/3}/\nu) + O(t^{-1/3}\log t)$, and hence the lower bound of (87) can be estimated as

$$\begin{aligned}
&\sqrt{\frac{b_2}{\log b_1 + 1/2\log(b_2/(\log b_1 - \log x)) - \log x}} \\
&\qquad = \frac{4t\ell_1(1+O(t^{-4/9}))}{(4t\ell_1/[\gamma c_2 t^{1/3}/\nu]_-^2) + \log(2\gamma c_2/\ell_1 s_1 \nu) - \log x}.
\end{aligned}$$

Substituting this lower bound on the right-hand side of (87), we get an upper bound for $y_0$ of the same kind, so that the quantity $\upsilon(l_1 + 1 > [y_0])$ equals

$$\begin{aligned}
&\exp\left\{-\nu\left[(1+O(t^{-4/9}))\right.\right. \\
&\qquad \left.\left.\times\sqrt{\frac{4t\ell_1}{(4t\ell_1/[\gamma c_2 t^{1/3}/\nu]_-^2) + \log(2\gamma c_2/\ell_1 s_1\nu) - \log x}}\right]\right\} \\
&\quad = \exp\left\{-\nu\left[\sqrt{\frac{4t\ell_1}{(4t\ell_1/[\gamma c_2 t^{1/3}/\nu]_-^2) + \log(2\gamma c_2/\ell_1 s_1\nu) - \log x}} + o(1)\right]\right\}.
\end{aligned} \tag{90}$$

Now, from the expansion $(1+y)^{-1/2} = 1 - \frac{1}{2}y + o(y^2)$ for small $y$, the choice of $s_1$ so that $\frac{2\gamma c_2}{\ell_1 s_1 \nu} = 1$, the fact that $n(t) = \frac{2cp}{1-p} e^{\nu[(\gamma c_2/\nu)t^{1/3}]_-} + O(1)$ and (90) combined with (89) for $\upsilon(Y_1(t) > x, t^{2/9} < l_1 < t^{2/5})$, we see that

$$\begin{aligned}
&n\upsilon(Y_1 > x, t^{2/9} < l_1 < t^{2/5}) \\
&\quad = 2c\frac{p}{1-p} \\
&\qquad \times \exp\left\{\nu\left[\frac{\gamma c_2}{\nu}t^{1/3}\right]_-\right. \\
&\qquad\qquad \left. - \nu\left[\left[\frac{\gamma c_2}{\nu}t^{1/3}\right]_-\left(1 + \frac{1}{2}\frac{[\gamma c_2 t^{1/3}/\nu]_-^2}{4t\ell_1}\log x\right) + o(1)\right]\right\} + o(1) \\
&\quad = 2c\frac{p}{1-p}\exp\left\{-\nu\left[\frac{1}{2}\frac{[\gamma c_2 t^{1/3}/\nu]_-^3}{4t\ell_1}\log x + o(1)\right]\right\} + o(1).
\end{aligned} \tag{91}$$

Now, note that $\frac{1}{2}\frac{[\gamma c_2 t^{1/3}/\nu]_-^2}{4t\ell_1} < \frac{1}{\nu}(\gamma\frac{3}{2})^3 = \frac{1}{\nu}a_1(\gamma)$ and $\lim_{t\to\infty}\frac{1}{2}\frac{[\gamma c_2 t^{1/3}/\nu]_-^2}{4t\ell_1} = \frac{1}{\nu}a_1(\gamma)$. Hence,

$$\lim_{t\to\infty} n\upsilon(Y_1 > x,\ t^{2/9} < l_1 < t^{2/5}) = \frac{2cp}{1-p}(1-p)^{-[(1/\nu)\log x^{a_1(\gamma)}]}.$$

Combining with (86), this proves (84).

6.4. *The Lévy spectral function in the annealed non-Gaussian case.* Now we compute the Lévy spectral function of the limiting law of part (ii) of Proposition 2. This time, we must verify the condition (81) with

$$Y_1(t) := \frac{1}{s_1 t^{1/3}} \exp\left\{\frac{4t\ell_1}{[\gamma c_2 t^{1/3}/\nu]_-^2}\right\}(X_1(t) - \upsilon(X_1(t))),$$

$L(t) := e^{\nu[(\gamma c_2/\nu)t^{1/3}]_-}$, $n(t) := 2r(t)+1$, $r(t) := [cL(t)\frac{p}{1-p}]$ and $c > 0$. Now, it is easy to check that

$$\lim_{t\to\infty} \frac{1}{c_2(1,p)t^{1/3}} \log \upsilon(X_1(t)) = -1.$$

Therefore, we have that

$$\begin{aligned}&\frac{1}{s_1 t^{1/3}} \exp\left\{\frac{4t\ell_1}{[\gamma c_2 t^{1/3}/\nu]_-^2}\right\}\upsilon(X_1(t))\\ &\qquad = \exp\{(a(\gamma) - \gamma - (1+\gamma))t^{1/3} + o(t^{1/3})\}.\end{aligned}$$

But, $(a(\gamma) - \gamma - (1+\gamma)) = \frac{2^2}{3^3}\frac{1}{\gamma^2} - (1+\gamma) < 0$, whenever $\gamma > \gamma_1 = 1/3$. It follows that for real $x$,

$$\upsilon(Y_1(t) > x) = \upsilon\left(\frac{1}{s_1 t^{1/3}} \exp\left\{\frac{4t\ell_1}{[\gamma c_2 t^{1/3}/\nu]_-^2}\right\}X_1(t) > x + o(e^{-t^{1/4}})\right).$$

It turns out that the term $o(e^{-t^{1/4}})$ is small enough to not affect the computation of the Lévy spectral function. In fact, essentially a repetition of the calculations of the previous subsection show that

$$\begin{aligned}\mathcal{L}(x) &= -\lim_{t\to\infty} n(t)\upsilon(Y_1(t) > x)\\ &= \begin{cases} 0, & \text{for } x < 0,\\ -\dfrac{2cp}{1-p}(1-p)^{-[(1/\nu)\log x^{a_1(\gamma)}]}, & \text{for } x > 0.\end{cases}\end{aligned}$$

6.5. *The truncated moments.* The purpose of this subsection is the verification of conditions (82) and (83) respectively showing that there is no spectral atom and identifying the constants $\beta_1$ and $\beta_2$ of the Lévy representation of the characteristic function of the limiting laws, in parts (i) and (ii) of Proposition 2. Namely we will prove that

$$\lim_{\tau\to+0}\lim_{t\to\infty} n(t)(\upsilon(Y(t)^2, Y(t)\le\tau) - \upsilon(Y(t), Y(t)<\tau)^2) = 0,$$

$$\beta_1 = \frac{2p}{(1-p)}\sum_{k=-\infty}^{\infty}\frac{(1-p)^k}{(1-p)^{k/a_1}+(1-p)^{-k/a_1}},$$

and that

$$\beta_2 = \frac{2p}{(1-p)}\sum_{k=-\infty}^{\infty}\frac{(1-p)^{k(1+2/a_1)}}{(1-p)^{k/a_1}+(1-p)^{-k/a_1}}.$$

Our first step is the following lemma.

LEMMA 11. *For every integer $k\ge 1$ it is true that*

$$\begin{aligned}
&\lim_{t\to\infty} n(t)\upsilon(Y^k(t), Y(t)\le\tau)\\
(92)\qquad &= 2c\left(\frac{p}{1-p}\right)^2 (1-p)^{(k-a_1)/a_1[(1/\nu)\log\tau^{a_1}]}\frac{1}{1-(1-p)^{1/a_1-1}}.
\end{aligned}$$

PROOF. Note that

$$\begin{aligned}
&n(t)\upsilon(Y^k(t), Y(t)\le\tau)\\
(93)\qquad &\sim \frac{2cp}{1-p}e^{\nu[(\gamma c_2/\nu)t^{1/3}]_-}\frac{\exp\{k4t\ell_1/[\gamma c_2 t^{1/3}/\nu]_-^2\}}{(s_1 t^{1/3})^k}\sum_{j=0}^{m} e^{-f_t(j)},
\end{aligned}$$

where

$$f_t(x) := \frac{4kt\ell_1}{(x+1)^2}(1+o_1(x)) + \nu x - k\log\frac{2(x+1)}{\ell_1} - \log(1+o_4(x,t))$$

and $m$ is the largest integer such that $\frac{e^{ac_2t^{1/3}}}{s_1t^{1/3}(2L+1)}\frac{2(m+1)}{\ell_1}e^{-(4t\ell_1/(m+1)^2)(1+o_1(m))}(1+o_4(m,t))\le\tau$. Let us now interchange the order in the summation of (93) so that $\sum_{j=0}^m e^{-f_t(j)} = \sum_{j=0}^m e^{-f_t(m-j)}$ and expand $f_t(m-j)$ around $m$,

$$f_t(m-j) = f_t(m) + \left(\frac{8tk\ell_1}{(m+1)^3}-\nu\right)j + O_1(j,m,t),$$

where

$$\begin{aligned}
O_1(j,m,t) &:= 4kt\ell_1\left(\frac{j^2}{(m+1)^4}\right)\left(\frac{1}{(1-j/(m+1))^2}+2\frac{1}{1-j/(m+1)}\right)\\
&\quad + \frac{o_1(m-j)4kt\ell_1}{(m+1-j)^2} + \log\left(1-\frac{j}{m+1}\right) - \log(1-o_4(m-j,t)).
\end{aligned}$$

Thus, we see that the right-hand side of (93) can be expressed as

$$\frac{2cp}{1-p} e^{\nu[(\gamma c_2/\nu)t^{1/3}]_-} \frac{\exp\{k4t\ell_1/[\gamma c_2 t^{1/3}/\nu]_-^2\}}{(s_1 t^{1/3})^k} \times e^{-f_t(m)} \sum_{j=0}^{m} e^{-((8tk\ell_1/(m+1)^3)-\nu)j+O_1(j,m,t)}.$$

But, by the definition of $m$, it is easy to see that

$$\lim_{t\to\infty} \frac{2cp}{1-p} e^{\nu[(\gamma c_2/\nu)t^{1/3}]_-} \frac{\exp\{k4t\ell_1/[\gamma c_2 t^{1/3}/\nu]_-^2\}}{(s_1 t^{1/3})^k} e^{-f_t(m)} = \frac{2cp^3}{(1-p)^3}(1-p)^{((k-a_1)/a_1)[(1/\nu)\log \tau^{a_1}]}.$$

Now, a straightforward computation shows that

$$\lim_{t\to\infty} \sum_{j=0}^{m} e^{-((8tk\ell_1)/((m+1)^3)-\nu)j+O_1(j,m,t)} = \frac{1}{1-(1-p)^{(1/a_1)-1}}. \qquad \square$$

Let us now compute $\beta_1$ via (83). Note that since $\int_0^\tau x\,d\mathcal{L}$ is well defined when $\gamma<\gamma_1$, we can make the decomposition $\int_{|x|<\tau} \frac{x^3}{1+x^2}\,d\mathcal{L} = \int_0^\tau x\,d\mathcal{L} - \int_0^\tau \frac{x}{1+x^2}\,d\mathcal{L}$. It follows that

$$\beta_1 = \lim_{t\to\infty} n(t)E_t(Z_\tau(t)) - \int_0^\tau x\,d\mathcal{L} + \int_0^\infty \frac{x}{1+x^2}\,d\mathcal{L}.$$

For $n\in\mathbb{Z}$ let $x_n := (1-p)^{n/a_1}$. Note that $\{x_n : n\in\mathbb{Z}\}$ are the discontinuities of $\mathcal{L}(x)$ and $[\frac{1}{\nu}\log x_n^{a_1}] = n$. We then have that

$$\begin{aligned}\int_0^\tau x\,d\mathcal{L} &= \sum_{k=m}^{\infty} x_k(\mathcal{L}(x_{k+1}) - \mathcal{L}(x_k)) \\ &= \frac{2cp}{1-p}\sum_{k=m}^{\infty}(1-p)^{k/a_1}((1-p)^{-k-1} - (1-p)^{-k}) \\ &= \lim_{t\to\infty} n(t)E_t(Z_\tau(t)),\end{aligned}$$

where $m = [\frac{1}{\nu}\log\tau^{a_1}]$ and we used Lemma 11. Hence,

$$\begin{aligned}\beta_1 = \int_0^\infty \frac{x}{1+x^2}\,d\mathcal{L} &= \sum_{n=-\infty}^{\infty} \frac{x_n}{1+x_n^2}(\mathcal{L}(x_{n+1}) - \mathcal{L}(x_n)) \\ &= \frac{2cp}{1-p}\sum_{k=-\infty}^{\infty} \frac{(1-p)^{k/a_1}}{1+(1-p)^{2k/a_1}}((1-p)^{-k-1} - (1-p)^{-k}).\end{aligned}$$

A similar computation enables us to compute $\beta_2$.

G. Ben Arous
Courant Institute of
Mathematical Sciences
New York University
New York, New York 10012
USA
e-mail: benarous@cims.nyu.edu

S. Molchanov
Department of Mathematics
University of North Carolina–Charlotte
376 Fretwell Bldg. 9201
University City Blvd.
Charlotte, North Carolina 28223
USA
e-mail: smolchan@math.uncc.edu

A. F. Ramírez
Facultad de Matemáticas
Pontificia Universidad Católica de Chile
Vicuña Mackenna 4860
Macul
Santiago 6904411
Chile
e-mail: aramirez@mat.puc.cl